\documentclass[12pt]{amsart}
\usepackage{amsmath,amsfonts,amssymb,amsthm, multicol}
\usepackage{anysize}
\marginsize{2cm}{2cm}{2cm}{2cm}
\usepackage{graphicx}
\usepackage{color}
\usepackage[pdftex]{hyperref}
\usepackage{enumerate}
\usepackage{caption}

\usepackage{mathtools}

\def\u{\mathfrak{u}}

\def\g{\mathfrak{g}}
\def\h{\mathfrak{h}}
\def\n{\mathfrak{n}}

\def\R{\mathbb{R}}

\def\Z{\mathbb{Z}}
\def\N{\mathbb{N}}

\def\e{\operatorname{e}}

\def\ad{\operatorname{ad}}
\def\tr{\operatorname{tr}}
\def\I{\operatorname{I}}
\def\alt{\raise1pt\hbox{$\bigwedge$}}

\def\vol{\operatorname{vol}}

\def\mid{\, \vert \,} 
\def\Span{\operatorname{span}}
\def\spec{\operatorname{spec}}

\theoremstyle{plain}
\newtheorem{theorem}{\bf Theorem}[section]
\newtheorem{corollary}[theorem]{\bf Corollary}
\newtheorem{proposition}[theorem]{\bf Proposition}
\newtheorem{lemma}[theorem]{\bf Lemma}

\theoremstyle{definition}
\newtheorem{definition}[theorem]{\bf Definition}
\newtheorem{example}[theorem]{\bf Example}

\theoremstyle{remark}
\newtheorem{remark}[theorem]{\bf Remark}

\title{Symplectic solvmanifolds not satisfying the hard-Lefschetz condition}

\author{Adrián Andrada}
\email{adrian.andrada@unc.edu.ar}
\author{Agustín Garrone}
\email{agustin.garrone@unc.edu.ar}

\date{}
\address{FAMAF, Universidad Nacional de C\'ordoba and CIEM-CONICET, Av. Medina Allende s/n, Ciudad Universitaria, X5000HUA C\'ordoba, Argentina}


\begin{document}

\begin{abstract}      
    For Lie groups $G$ of the form $G = \R^k \ltimes_{\phi} \R^m$, with $k + m$ even, a result of H. Kasuya shows that if the action $\phi:\R^k \to \mathrm{Aut}(\R^m)$ is semisimple then any symplectic solvmanifold $(\Gamma \backslash G, \omega)$ satisfies the hard-Lefschetz condition for any symplectic form. In this article, we prove the converse in the case $k = 1$ and $G$ completely solvable: no symplectic form on such a solvmanifold satisfies the hard-Lefschetz condition if $\phi$ is not semisimple; moreover, we show that the failure occurs either at degree $1$ or at degree $2$ in cohomology, depending on the spectrum of the differential of the action $\phi$. This result is achieved through a detailed analysis of the cohomology groups $H^1(\g)$, $H^2(\g)$, $H^{2n-2}(\g)$, $H^{2n-1}(\g)$ of the Lie algebra $\g$ of such Lie groups. Among other things, this analysis yields useful representatives for each cohomology class corresponding to any symplectic form on $\g$, allowing the most delicate cases to be reduced to a straightforward computation. We also construct lattices for many of the Lie groups under consideration, thereby exhibiting examples of symplectic solvmanifolds of completely solvable Lie groups failing to have the hard-Lefschetz property for any symplectic form. 
\end{abstract}

\maketitle

\section{Introduction} \label{section: introduction}

\indent The hard-Lefschetz condition is a stringent cohomological property shared by all compact Kähler manifolds, and can be loosely described as the possibility of developing an analogue of Hodge theory for symplectic manifolds exhibiting said property. We elaborate on these concepts and in its importance in Section \ref{section: HLC}. For the time being, it suffices to mention that the hard-Lefschetz condition guides the understanding of Kähler-like manifolds that are not Kähler, at least from a cohomological point of view. 

\indent On a compact symplectic manifold $(M,\omega)$ of dimension $2n$, whose algebra of differential forms is denoted by $\Omega(M)$, the hard-Lefschetz condition can be stated as the fact that the \textit{Lefschetz operators} $L_m:\Omega^{n-m}(M) \to \Omega^{n+m}(M)$, given by $L_m(\alpha) := \omega^m \wedge \alpha$, induce isomorphisms $H_{dR}^{n-m}(M) \to H_{dR}^{n+m}(M)$ on the de Rham cohomology for all $0 \leq m \leq n$. In this article we are particularly interested in the situation when $M$ is a \textit{solvmanifold}, which is just a quotient $\Gamma \backslash G$ of some simply connected solvable Lie group $G$ by a discrete co-compact subgroup $\Gamma$ of $G$, called a \textit{lattice}. A short review of solvmanifolds is given in Section \ref{section: solvmanifolds}. The existence of lattices is difficult to establish, and not many more necessary conditions are known besides the fact that $G$ should be unimodular (due to a result by Milnor that we collect in Proposition \ref{prop: Milnor}). When $G$ is a completely solvable Lie group, there is an isomorphism of cohomology rings $H^*(\g) \cong H^*_{dR}(\Gamma\backslash G)$ induced by the natural inclusion $\alt^* \g^* \hookrightarrow \Omega^*(\Gamma\backslash G)$, allowing the study of cohomological properties at the Lie algebra level \cite{Hattori}. We focus on \textit{almost abelian} Lie groups $G$, defined as the simply connected Lie groups corresponding to a Lie algebra $\g$ that has a codimension-one abelian ideal. These Lie algebras can be described alternatively as semidirect products $\R \ltimes_A \R^m$, which we denote by $\g_A$. In this case, the corresponding simply connected Lie group $G_A$ of $\g_A$ is also a semidirect product $\R \ltimes_{\phi} \R^m$, with $\phi(t) = \exp(tA)$. Preliminaries on this topic are included in Section \ref{section: almost abelian lie algebras}. We mention here only that the semisimplicity of the action $\phi$ is equivalent to the fact that $A$ is a semisimple matrix. 

\indent The main result of this article is as follows. Remark \ref{obs: thank you, hattori} below may be relevant.   
\begin{theorem} \label{thm: main result} 
    Let $\Gamma \backslash G$ be a symplectic completely solvable solvmanifold with $G = \R \ltimes_{\phi} \R^m$. If the action $\phi$ is not semisimple then no symplectic form on $\Gamma \backslash G$ satisfies the hard-Lefschetz condition. 
\end{theorem}


\indent Moreover, as stated in Theorems \ref{thm: main result 1} and \ref{thm: main result 2}, we are able to show which is the first Lefschetz operator that fails to be an isomorphism. Theorem \ref{thm: main result} is built conceptually upon Theorems \ref{thm: Benson y Gordon} and \ref{thm: Kasuya}, due to Benson-Gordon and Kasuya, respectively. Theorem \ref{thm: Benson y Gordon} is actually used as a part of the proof, and together with Theorem \ref{thm: Kasuya} one gets the more general statement. 
\begin{theorem}
    Let $\Gamma \backslash G$ be a symplectic completely solvable solvmanifold with $G = \R \ltimes_{\phi} \R^m$. If the action $\phi$ is semisimple then any symplectic form on $\Gamma \backslash G$ satisfies the hard-Lefschetz condition, whereas if the action $\phi$ is not semisimple then no symplectic form on $\Gamma \backslash G$ satisfies the hard-Lefschetz condition.
\end{theorem}
\indent Note that there is no in-between when it comes to this kind of solvmanifolds: Either none or all of the symplectic forms are hard-Lefschetz. Having in mind \cite[Theorem 1.1]{CM} and our own result Theorem \ref{thm: splitting of symplectic forms}, discussed in Remark \ref{obs: cohomological uniqueness?}, this may suggest some sort of ``essential uniqueness"\! for symplectic forms in cohomology. Having no more evidence to support this, we only wish to pose a modest question, of which Theorem \ref{thm: main result} constitutes the answer for the case $k = 1$, in the particular case of completely solvable groups. \\

\noindent \textbf{Question.} Assume that the solvmanifold $\Gamma \backslash G$ with $G = \R^k \ltimes_{\phi} \R^m$ admits symplectic forms. Is it true that if the action $\phi:\R^k \to \mathrm{Aut}(\R^m)$ is \underline{not} semisimple then $(\Gamma\backslash G,\omega)$ does not satisfy the hard-Lefschetz condition, regardless of the choice of the symplectic form? \\
 
\indent The proof of Theorem \ref{thm: main result} takes about all Section \ref{section: main section}. Not much can be said without the special notation introduced in Section \ref{section: notation}. The key steps are Proposition \ref{prop: splitting of symplectic forms} and Theorem \ref{thm: splitting of symplectic forms}, the first being a consequence of the careful study done in Section \ref{section: the space of closed 2 forms}, arguably the bulk of our efforts, and the second being mainly a trick. They essentially describe how any symplectic form on $\g_A$, under mild conditions, splits first as a sum of symplectic forms on certain subspaces we call \textit{double generalized eigenspaces}, and then as a sum of closed $2$-forms in further subspaces reflecting the Jordan form of $A$. All our proofs are elementary and, in many cases, outright computation-oriented. Lastly, in Section \ref{Section: lattices} we briefly discuss the existence of lattices for some of the Lie groups we consider, thereby exhibiting examples of symplectic solvmanifolds not satisfying the hard-Lefschetz condition. 

\section{Preliminaries} \label{section: preliminaries} 

\subsection{The hard-Lefschetz condition.} \label{section: HLC}

\indent Let $(M,\omega)$ be a compact symplectic manifold of dimension $2n$, and denote by $\Omega(M)$ the algebra of differential forms in $M$. There is a natural volume form on $M$ given by $\vol_{\omega} = \frac{\omega^n}{n!}$, where $\omega^n := \omega \wedge \cdots \wedge \omega$ ($n$ times). The symplectic form also induces a skew-symmetric map
\begin{align*}
    \omega^{-1}: \Omega(M) \times \Omega(M) \to \mathbb{R}, \quad \omega^{-1}(\theta, \eta) :=\eta( B^{-1}(\theta)),
\end{align*}
\noindent where $B^{-1}:\Omega^1(M) \to \mathfrak{X}(M)$ is the inverse of
\begin{align*}
    B:\mathfrak{X}(M) \to  \Omega^1(M), \quad B(X)(Y) := \omega(X,Y)
\end{align*}
\noindent whose existence is due to the non-degeneracy of $\omega$. We can extend $\omega^{-1}$ to a bilinear form on $\Omega(M)$ given by the $C^{\infty}(M)$-linear extension of 
\begin{align*}
    \omega^{-1}( \theta_1 \wedge \cdots \wedge \theta_k, \eta_1 \wedge \cdots \wedge \eta_l ) = 
    \begin{cases}
        \det( \left[ \omega^{-1}(\theta_i, \eta_j) \right]_{1 \leq i,j \leq k} ), & k = l,\\
        0 & k \neq l,
    \end{cases}
\end{align*}
\noindent and $\theta_i$'s and $\eta_j$'s are $1$-forms. With this we can build up a rudimentary Hodge theory on $(M,\omega)$, as it is possible to define a \textit{symplectic star operator} $\star_{\omega}:\Omega(M) \to \Omega(M)$, 
\begin{align*}
    \alpha \wedge \star_{\omega} \beta = \omega^{-1}(\alpha,\beta) \vol_{\omega}, 
\end{align*}
\noindent and a \textit{symplectic codifferential} $d^c:\Omega(M) \to \Omega(M)$, 
\begin{align*}
   d^c\vert_{\Omega^k(M)} = (-1)^{k+1} \star_{\omega} \circ \, d \circ \star_{\omega}. 
\end{align*}
\noindent Certainly, $\star_{\omega}$ is a linear involution such that $\star_{\omega}(\Omega^k(M)) \subseteq \Omega^{2n - k}(M)$ for all $0 \leq k \leq 2n$, and $d^c$ is a differential (i.e. a linear map with $d^c \circ d^c = 0$) that skew-commutes with $d$ and satisfies $d^c(\Omega^k(M)) \subseteq \Omega^{k-1}(M)$ for all $0 \leq k \leq 2n$. It is straightforward to verify that
\begin{align} \label{eq: ddc lemma}
    \star_\omega(\operatorname{Ker} d^c \cap \mathrm{Im} \; d) = \operatorname{Ker} d \cap \mathrm{Im} \; d^c. 
\end{align}
\noindent Since $d d^c = - d^c d$, the natural candidate for a symplectic Laplacian $\Delta_{\omega} = d d^c + d^c d$ is identically zero. We thus say that a differential form $\alpha \in \Omega(M)$ is \textit{$\omega$-harmonic} if it is $d$-closed and $d^c$-closed; that is, if $d \alpha = 0$ and $d^c \alpha = 0$.  


\indent As in the case of compact K\"ahler manifolds, the exterior algebra $\Omega(M)$ of $(M, \omega)$ supports a representation of the Lie algebra $\mathfrak{sl}(2, \mathbb{R})$ given by the triple of operators $L$, $\Lambda$, $H:\Omega(M) \to \Omega(M)$ defined as
\begin{align*}
    L(\alpha) := \omega \wedge \alpha, \quad \Lambda(\alpha) := \star_{\omega} \circ L \circ \star_{\omega} \, \alpha, \quad H := [L, \Lambda].
\end{align*}
\noindent Note that, for all $0 \leq k \leq n$, we have 
\begin{align*}
    L(\Omega^k(M)) \subseteq \Omega^{k+2}(M), \quad \Lambda(\Omega^k(M)) \subseteq \Omega^{k-2}(M), \quad H(\Omega^k(M)) \subseteq \Omega^k(M).
\end{align*}
\noindent The repeated applications of the $L$-operator,
\begin{align*}
    L_k:\Omega^{n-k}(M) \to \Omega^{n+k}(M), \quad L_k(\alpha) = \omega^k \wedge \alpha,
\end{align*}
\noindent induce operators in de Rham cohomology, which we call by the same name in a slight abuse of language as it is customary, by
\begin{align*}
    L_k: H_{dR}^{n-k}(M) \to H_{dR}^{n+k}(M), \quad L_k([\alpha]) = [\omega^k \wedge \alpha]. 
\end{align*}
\noindent These operators are known as \textit{Lefschetz operators}, and are crucial to us as the following result summarizes.
\begin{theorem} \cite[Theorem 3.1]{TT} \label{thm: HLC equivalences}
    The following properties are equivalent on any compact symplectic manifold $(M, \omega)$.
\begin{enumerate} [\rm (i)]
    \item The hard-Lefschetz condition: The Lefschetz operators $L_k$ are linear isomorphisms for all $0 \leq k \leq \frac{\dim M}{2}$.  
    \item Hodge theorem: Any de Rham cohomology class on $M$ has a $\omega$-harmonic representative.
    \item The $dd^c$-lemma: Any $d^c$-closed and $d$-exact form is $d d^c$-exact; equivalently, any $d$-closed and $d^c$-exact form is $dd^c$-exact. In symbols, the following equivalent\footnote{This equivalence holds because of equation \eqref{eq: ddc lemma}.} identities hold:
\begin{align*}
    \operatorname{Ker} d^c \cap \mathrm{Im} \, d = \mathrm{Im} \, d d^c \quad \Longleftrightarrow \quad \operatorname{Ker} d \cap \mathrm{Im} \, d^c = \mathrm{Im} \, d d^c.
\end{align*}
\end{enumerate} 
\end{theorem}
\indent As indicated in \cite[Theorem 3.1]{TT}, Theorem \ref{thm: HLC equivalences} is the culmination of works by several authors spanning several years of research. In fact, Brylinski \cite[Conjecture 2.2.7]{Brylinski} conjectured that (i) and (ii) are equivalent. A few years later, Mathieu \cite[Corollary 2]{Mathieu} and Yau \cite[Theorem 0.1]{Yan} independently proved this conjecture. Later, Merkulov \cite[Proposition 1.4]{Merkulov} established the equivalence of (ii) and (iii), a revised version of which was published by Cavalcanti \cite[Theorem 5.4]{Cavalcanti 1}. 
\begin{remark}
    Some authors define the \textit{Lefschetz condition} (without the ``hard") as the requirement that the Lefschetz operator $L_{n-1}:H_{dR}^1 \to H_{dR}^{2n-1}$ at degree $1$ is a linear isomorphism. See, for example, \cite{Yamada 1}, where some examples of symplectic manifolds satisfying the Lefschetz but not the hard-Lefschetz condition are shown. Lefschetz \textit{solvmanifolds} (see Section \ref{section: solvmanifolds}) are characterized in \cite[Theorem 2]{BG2}. 
\end{remark}
\indent Being the main inspiration behind the concept, certainly all K\"ahler manifolds satisfy the hard-Lefschetz condition. On the other hand, many compact symplectic manifolds exhibit that condition but fail to be K\"ahler. To our knowledge, the first historically significant examples of such a manifold were produced either by Cavalcanti \cite[Section 4]{Cavalcanti 2} or by de Bartomoleis and Tomassini \cite{BT}, both being surprisingly recent. Many other examples are known in the literature, including those described in our previous work \cite{AG}. Remarkably, the fact that a compact manifold admits K\"ahler structures does not imply that all other symplectic forms satisfy the hard-Lefschetz condition, as the following result shows. 
\begin{proposition} \cite[Theorem 1.3]{Cho} \label{prop: Cho}
    There exists a simply connected compact K\"ahler manifold $(X,\omega)$ that admits a symplectic form $\sigma$ such that $(X, \sigma)$ does not satisfy the hard-Lefschetz condition. Moreover, $X$ has real dimension 6, its odd Betti numbers are all zero, and $\sigma$ is deformation-equivalent to $\omega$.
\end{proposition}

\indent We mention briefly that Theorem \ref{thm: HLC equivalences} links the hard-Lefschetz condition to the possibility of studying compact symplectic manifolds via several cohomological theories, some of which are inspired by counterparts in almost complex geometry: we are talking about \textit{Aeppli} and \textit{Bott-Chern} cohomologies. See, for instance, \cite{AK}, as well as the references cited therein. In addition to this and the aforementioned relation to the feasibility of developing an analog of Hodge theory for compact symplectic manifolds, it is known that symplectic manifolds satisfying the hard-Lefschetz condition have Betti numbers that possess a number of properties commonly referred to as ``the Kähler package": 
\begin{itemize}
    \item The first property in the package is that odd Betti numbers are even. This is a consequence of the maps $B:H^{2i+1}_{dR}(M) \times H^{2i+1}_{dR}(M) \to \mathbb{R}$ given by $B([\alpha], [\beta]) := \star_{\omega} L_{n - 2i - 1}([\alpha \wedge \beta])$ being non-degenerate skew-symmetric forms on $H^{2i+1}_{dR}(M)$ for all $0 \leq i \le n-1$. 
    \item The second property in the package is that Betti numbers are \textit{unimodal}, meaning that they form an increasing sequence up to degree $\frac{\dim M}{2}$, and decreasing from then on. This is a consequence of the Poincaré duality, as well as the maps $C:H^i_{dR}(M) \to H^{i+2}_{dR}(M)$ given by $C([\alpha]) = [\omega \wedge \alpha]$ being injective for all $0 \leq i \leq n-2$. 
    \item The third property in the package is that even Betti numbers are strictly positive. This is actually not a consequence of the hard-Lefschetz condition, but rather of compactness and the Stokes theorem (as otherwise $\frac{\omega^n}{n!}$ would not be a volume form on $M$). 
\end{itemize}

\indent We close this subsection by stating and proving the following elementary result for future reference. We informally refer to it as the ``propagation of the non-Lefschetz condition".
\begin{lemma} \label{lemma: propagacion}
    Let $(M, \omega)$ be a direct product of compact symplectic manifolds $(M_1, \omega_1)$ and $(M_2, \omega_2)$. If $(M_1, \omega_1)$ is not hard-Lefschetz then neither is $(M, \omega)$. Moreover, if the Lefschetz condition fails in $(M_1, \omega_1)$ at degree $k$ \footnote{If $(N^{2n},\rho)$ is a symplectic manifold, we say that the hard-Lefschetz condition fails at degree $k$ if the Lefschetz operator $L_{n - k}:H^k(N) \to H^{2n - k}(N)$ has nontrivial kernel.} then it also fails for $(M, \omega)$ at degree $k$.  
\end{lemma}
\begin{proof}
    \noindent Set $n_1,$ $n_2 \in \N$ such that $\dim M_1 = 2 n_1$ and $\dim M_2 = 2n_2$, and define $n := n_1 + n_2$. Assume that there exists $[\gamma] \in H^k_{dR}(M_1)$ such that $[\omega_1^{n_1 - k} \wedge \gamma] = 0$. Note that there is an (injective) inclusion $\iota: H^k_{dR}(M_1) \to H^k_{dR}(M)$ induced by the canonical projection $\pi:M \to M_1$. Denote $\iota([\gamma])$ simply as $[\gamma]$; in similar fashion, write $\omega$ just as $\omega_1 + \omega_2$. By the binomial theorem, and given that $\omega_1^{n_1}$ and $\omega_2^{n_2}$ are top forms on $M_1$ and $M_2$ respectively, we have that
\begin{align*}
    \omega^{n - k} = (\omega_1 + \omega_2)^{n_1 + n_2 - k} = \sum_{r = 0}^{n-k} \binom{n-k}{r} \omega_1^r \wedge \omega_2^{n-k-r} = \sum_{r = n_1 - k}^{ \min\{n_1, n - k\}} \binom{n-k}{r} \omega_1^r \wedge \omega_2^{n-k-r}.
\end{align*}
    \noindent Therefore,
\begin{align*}
    \omega^{n - k} \wedge \gamma &= \left( \sum_{r = n_1 - k}^{\min \{n_1, n - k\}} \binom{n-k}{r} \omega_1^r \wedge \omega_2^{n-k-r} \right) \wedge \gamma \\
    &= \left( \sum_{r = 0}^{\min\{k, n_2\}} \binom{n-k}{n_1 - k + r} \omega_1^{n_1 - k + r} \wedge \omega_2^{n_2 - r} \right) \wedge \left( \omega_1^{n_1 - k} \wedge \gamma \right). 
\end{align*}
    \noindent As $\sum_{r = 0}^{\min \{k, n_2\}} \binom{n-k}{n_1 - k + r} \omega_1^{n_1 - k + r} \wedge \omega_2^{n_2 - r}$ is closed and $\omega_1^{n_1 - k} \wedge \gamma$ is exact, $\omega^{n - k} \wedge \gamma = 0$ is also exact. 
\end{proof}
\indent An adaptation of Lemma \ref{lemma: propagacion} to the Lie algebra setting, more suited to our purposes, is stated and proved in Section \ref{section: almost abelian lie algebras}; see Lemma \ref{lemma: propagacion algebraica}. 

\indent The rest of the article deals with the hard-Lefschetz condition on a specific kind of compact manifolds called \textit{solvmanifolds}. A review of them is in order. 

\subsection{Solvmanifolds.} \label{section: solvmanifolds}

Throughout the section, let $G$ denote a connected real Lie group with Lie algebra $\g$. 
\begin{definition}
    A \textit{solvmanifold} is a compact quotient $\Gamma \backslash G$, where $G$ is simply connected and solvable, and $\Gamma$ is a discrete subgroup of $G$. Such a co-compact discrete subgroup $\Gamma$ is called a \textit{lattice} of $G$.
\end{definition}
\indent Two special classes of solvmanifolds relevant in what follows are those arising respectively from nilpotent groups $G$, in which case $\Gamma \backslash G$ is called a \textit{nilmanifold}, and completely solvable groups $G$. Recall that a connected solvable Lie group $G$ is \textit{completely solvable} if all the adjoint operators $\ad_x: \g \to \g$, with $x\in \g$, have only real eigenvalues. 

\indent As every connected and simply connected solvable Lie group is diffeomorphic to $\R^n$, the usual argument involving the long exact sequence associated to a fibration implies that every solvmanifold $\Gamma \backslash G$ is aspherical, meaning that $\pi_n(\Gamma\backslash G) = 0$ for all $n > 1$, as well as $\pi_1(\Gamma \backslash G) = \Gamma$. Moreover, a classical result due to Mostow shows that solvmanifolds are \textit{rigid}, meaning that they are determined up to diffeomorphism by their fundamental groups. 
\begin{proposition} \cite[Theorem A]{Mostow} \label{prop: Mostow}
    If $\Gamma_1$ and $\Gamma_2$ are lattices in simply connected solvable Lie groups $G_1$ and $G_2$ respectively, and $\Gamma_1$ is isomorphic to $\Gamma_2$, then $\Gamma_1 \backslash G_1$ is diffeomorphic to $\Gamma_2 \backslash G_2$. 
\end{proposition} When $G$ is completely solvable, the rigidity result in Proposition \ref{prop: Mostow} can be strengthened. 
\begin{proposition} \cite[Theorem 5]{Saito} \label{prop: Saito}
   If $\Gamma_1$ and $\Gamma_2$ are lattices in simply connected completely solvable Lie groups $G_1$ and $G_2$ respectively, then every isomorphism $f: \Gamma_1 \to \Gamma_2$ extends uniquely to a Lie group isomorphism $F: G_1 \to G_2$. In particular, $\Gamma_1 \backslash G_1$ is diffeomorphic to $\Gamma_2 \backslash G_2$.
\end{proposition}

\indent The analogous results of Propositions \ref{prop: Mostow} and \ref{prop: Saito} for nilmanifolds were proven in \cite{Malcev} a few years earlier. 

\indent Solvmanifolds of completely solvable Lie groups have a very useful property concerning their de Rham cohomology, proved by Hattori, which is why they are relevant to our present discussion.
\begin{theorem} \cite[Theorem 4.1]{Hattori} \label{thm: Hattori}
    If $\Gamma \backslash G$ is a solvmanifold of a completely solvable Lie group $G$ then the natural inclusion $\alt^* \g^* \hookrightarrow \Omega^*(\Gamma\backslash G)$ induces an isomorphism $H^*(\g) \cong H^*_{dR}(\Gamma\backslash G)$. 
\end{theorem}
\indent That is, the de Rham cohomology of a solvmanifold of a completely solvable Lie group can be computed in terms of left invariant forms, and in particular its de Rham cohomology does not depend on the lattice $\Gamma$.
\begin{remark} \label{obs: first betti number}
    In the general case, the natural inclusion $\alt^* \g^* \hookrightarrow \Omega^*(\Gamma\backslash G)$ induces just an injective map $H^*(\g) \hookrightarrow H^*_{dR}(\Gamma\backslash G)$. Coupling this with the well-known facts that $\g \neq [\g,\g]$ for solvable $\g$ and that $H^1(\g) \cong \g/[\g,\g]$, we get that $b_1(\Gamma\backslash G)\geq b_1(\g) \geq 1$.
\end{remark} 
\indent The analogous result of Theorem \ref{thm: Hattori} for nilmanifolds was proven in \cite[Theorem 1]{Nomizu} a few years earlier. Results such as Proposition \ref{prop: Mostow} and Theorem \ref{thm: Hattori} seem to suggest that nilmanifolds are fitting testing grounds for conjectures concerning completely solvable solvmanifolds. Nonetheless, it is often the case that results valid in the class of nilmanifolds fail to hold for general solvmanifolds. In that vein, recall the following classical result due to Malcev. 
\begin{theorem} \cite{Malcev} \label{thm: Malcev} 
    A simply connected nilpotent Lie group $G$ admits lattices if and only if there is a basis of $\g$ with rational structure coefficients.  
\end{theorem}
\indent In contrast to Theorem \ref{thm: Malcev}, there appears to be no criterion for the existence of lattices in the class of general solvable Lie groups, not even completely solvable ones. The most we can hope for the goals of this article is a necessary criterion proved by Milnor. 

\begin{proposition} \cite[Lemma 6.2]{Milnor} \label{prop: Milnor}
    If a Lie group $G$ admits lattices then it is unimodular\footnote{An equivalent notion to unimodularity on Lie groups is that the Haar measure on $G$ is left and right invariant, a fact proved in \cite[Lemma 6.3]{Milnor}. }; i.e., $\tr(\ad_x) = 0$ for all $x \in \g$. 
\end{proposition}
\indent In this article we restrict to a specific class of completely solvable Lie groups for which there is a useful criterion to determine the existence of lattices: They are known as \textit{almost abelian Lie groups}, and can be characterized as semidirect products $G = \R \ltimes_{\phi} \R^m$ of abelian Lie groups. The criterion is as follows.
\begin{proposition} \cite[Corollary 4.5]{Bock} \label{prop: useful criterion for lattices}
    Let $G=\R\ltimes_\phi\R^m$ be a unimodular almost abelian Lie group. Then $G$ admits a lattice if and only if there exists $t_0\neq 0$ such that $\phi(t_0)$ is conjugate to a matrix in $\operatorname{SL}(m,\Z)$. In this situation, a lattice is given by $\Gamma = t_0 \Z \ltimes P\mathbb Z^m$, where $P\in \operatorname{GL}(m,\R)$ satisfies $P^{-1}\phi(t_0)P\in \operatorname{SL}(m,\Z)$. 
\end{proposition}

\indent We review that class of almost abelian Lie groups in Section \ref{section: almost abelian lie algebras}. While not particularly relevant in what follows, it is conceptually important to note that in any fixed dimension only countably many non-isomorphic simply
connected Lie groups admit lattices, according to Milovanov \cite[Theorem 4]{Milovanov} (for the solvable case) and Witte \cite[Proposition 8.7]{Witte} (for the general case).

\indent Another contrasting point between nilmanifolds and general solvmanifolds concerns the validity of the hard-Lefschetz condition.  
\begin{theorem} \cite[proof of Theorem A]{BG1} \label{thm: Benson y Gordon}
    A symplectic nilmanifold is hard-Lefschetz if and only if it is diffeomorphic to a torus, irrespective of the choice of the symplectic form. Moreover, the failure of the hard-Lefschetz condition for symplectic non-toral nilmanifolds occurs at degree $1$. 
\end{theorem}

\begin{remark} \label{obs: thank you, hattori}
    Usually, geometric structures in the context of solvmanifolds are assumed to be \textit{invariant}, meaning that they are induced from the corresponding left-invariant geometric structures defined on $G$ via the canonical projection map $G \to \Gamma \backslash G$ by imposing that it locally preserves the said structure. Certainly, one can do exactly this in the case of symplectic structures. A remarkable consequence of Hattori's Theorem \ref{thm: Hattori} is that there is no loss of generality in restricting to invariant symplectic forms when considering the cohomological properties of solvmanifolds associated to completely solvable Lie groups, as any symplectic form on such a solvmanifold is cohomologous to an invariant symplectic form.
\end{remark}

\indent In stark contrast to Theorem \ref{thm: Benson y Gordon}, the characterization of hard-Lefschetz symplectic solvmanifolds has proven to be challenging, and is currently far from complete even in the completely solvable case. An important result in this regard, closely related to the almost abelian Lie groups we are interested in, is the following. 
\begin{theorem} \cite[Corollary 1.5]{Kasuya} \label{thm: Kasuya}
    A symplectic solvmanifold $\Gamma \backslash G$ with $G = \R^k \ltimes_{\phi} \R^m$ and semisimple action $\phi\colon \R^k \to \mathrm{Aut}(\R^m)$ is hard-Lefschetz with respect to any symplectic form. 
\end{theorem}
\indent Notice the remarkable fact that this holds for every symplectic form, against the unsettling behavior described in Proposition \ref{prop: Cho} above. The example in this proposition is not a solvmanifold, as its odd Betti numbers all vanish, whereas solvmanifolds have a positive first Betti number as noted in Remark \ref{obs: first betti number}.

\indent In a previous work, we have explicitly determined the cohomology of two large families of solvmanifolds with $G = \R \ltimes_{\phi} \R^m$ and semisimple action $\phi$ by combinatorial means, establishing the same result as in Theorem \ref{thm: Kasuya} by direct computation, and exhibited explicit examples of lattices of such a group (see \cite{AG}). The articles \cite{Sawai} and \cite{Yamada 1} are concerned with similar (although different) examples; the same authors, in \cite{SY} and \cite{Yamada 2}, have also constructed lattices associated to Lie groups appearing in two of the examples given in \cite{BG2}, where a characterization of solvmanifolds of completely solvable Lie groups satisfying the $1$-Lefschetz property (i.e., those where only the map $L_{n-1}:H^1 \to H^{2n-1}$ is an isomorphism) are characterized. All these references are concerned with completely solvable Lie groups exclusively. For the non completely solvable case, the scarcity of examples is even more prominent. In this regard, we point out the interesting examples described in \cite{BT} and \cite{LT}.

\indent Combining Theorems \ref{thm: Benson y Gordon} and \ref{thm: Kasuya}, a natural question arises. Assume that the solvmanifold $\Gamma \backslash G$ with $G = \R^k \ltimes_{\phi} \R^m$ admits symplectic forms. Is it true that if $\phi:\R^k \to \mathrm{Aut}(\R^m)$ is \underline{not} semisimple then $(\Gamma\backslash G,\omega)$ does not satisfy the hard-Lefschetz condition, regardless of the choice of the symplectic form?

\indent In this article we answer the question for the affirmative in the case $k = 1$ under the additional hypothesis of complete solvability. We point out that in \cite{Kasuya-JDG} it is proven that the semisimplicity of the action $\phi$ above is in fact equivalent to the hard-Lefschetz condition for the cohomology of the solvmanifold with local systems; this result by its own is not enough to answer our question.  

\section{Almost abelian Lie algebras} \label{section: almost abelian lie algebras}

A finite dimensional real Lie algebra $\g_A$ is said to be \textit{almost abelian} if it has a codimension-one abelian ideal $\u$, from which it follows that it is a semidirect product $\g_A = \R \ltimes \u=\R\ltimes_A \R^m$ for some matrix $A$ of the appropriate size, after identifying $\u$ with $\R^m$ via a choice of basis: With this identification, $\g_A$ is the vector space $\R f_1 \oplus \u$ equipped with Lie bracket $[\cdot, \cdot]: \g_A \times \g_A \to \g_A$ given by 
\begin{align} \label{eq: bracket relations}
    [\u, \u] = 0 \quad \text{and} \quad \text{$[f_1, u] = Au$ for all $u \in \u$}.  
\end{align} 
\indent The corresponding simply connected Lie group $G_A$ of $\g_A$ is also a semidirect product $G_A = \R \ltimes_{\phi} \R^m$ with $\phi(t)=\exp(t A)$, and it is called \textit{almost abelian} as well. The matrix $A$ is called the \textit{structure matrix} of $\g_A$. It is clear from the definitions that the action $\phi$ is semisimple if and only if $A$ is a semisimple matrix. It is known that $\g_A$ is isomorphic to $\g_{A'}$ if and only if $cA$ and $A'$ are conjugate for some scalar $c\neq 0$; we refer to \cite[Proposition 2.9]{Freibert} for a proof. 

\indent With these definitions, it is clear that $\g_A$ is a $2$-step solvable Lie algebra; also, it is unimodular if and only if $\tr(A) = 0$, it is completely solvable if and only if $A$ has only real eigenvalues,  and it is nilpotent if and only if $A$ is nilpotent (moreover, $\u$ is the nilradical of $\g_A$ except when $A$ is nilpotent). The discussion in Section \ref{section: preliminaries} illustrates the importance of keeping track of these properties: Unimodularity is a necessary condition for lattices of $G_A$ to exist according to Proposition \ref{prop: Milnor}, complete solvability is necessary for Theorem \ref{thm: Hattori} to hold, and nilpotency is to be avoided if we intend to deal with cases not covered in Theorem \ref{thm: Benson y Gordon}. 

\indent Recall that the standard cochain complex $(\alt^\ast \g^\ast, d)$ in any Lie algebra $\g$ with bracket $[\cdot, \cdot]$ has exterior derivative $d$ completely determined by the conditions
\begin{gather*}
    d\theta(x,y) = - \theta([x,y]), \quad \text{for all $\theta \in \g^\ast$ and any $x$, $y \in \g$}, \\
    d(\alpha \wedge \beta) = d \alpha \wedge \beta + (-1)^{k} \alpha \wedge d \beta, \quad \text{for all } \alpha\in\alt^k\g^\ast, \, \beta \in \alt^\ast \g^\ast. \notag
\end{gather*}
\noindent Note in particular that, barring inaccuracies, the exterior derivative is minus the transpose of the Lie bracket. We mention in passing that an explicit formula for the exterior derivative of a generic $p$-form is possible, but we refrain from stating it in full generality since we are only concerned with the case $p = 2$: if $\eta \in \alt^2\g^\ast$ is a $2$-form on $\g$ then 
\begin{align*}
    - d\eta(x,y,z) = \eta([x,y],z) + \eta([y,z],x) + \eta([z,x],y) \text{ for all $x$, $y$, $z \in \g$}.  
\end{align*} 
\noindent A symplectic form on $\g$ is therefore a non-degenerate $2$-form $\omega \in \alt^2\g^\ast$ such that $d\omega = 0$, or equivalently
\begin{align}\label{eq: d2}
    \omega([x,y],z) + \omega([y,z],x) + \omega([z,x],y) = 0 \text{ for all $x$, $y$, $z \in \g$}. 
\end{align}
\indent Recall the following result (see also \cite[Proposition 4.1 and Remark 4.2]{LW}). 
\begin{proposition} \label{prop: LW} \cite[Proposition 4.2]{etal-nilp}
    Let $\omega \in \alt^2 \g_A^*$ a symplectic form on $\mathfrak{g}_A$. 
\begin{enumerate} [\rm (i)]
    \item There exists $f_2 \in \u$ such that $[f_1,f_2]= a f_2$ for some $a\in\R$ and $\omega(f_1,f_2)=1$.
    \item Let $\u_0$ be the $\omega$-orthogonal complement of $\{f_1, f_2\}$ in $\g_A$. Then $\u_0\subset \u$ and there exists a non-degenerate 2-form $\omega_0 \in \alt^2 \u_0^*$ such that $\omega = f^1 \wedge f^2 + \omega_0$. 
    \item According to the decomposition $\u=\R f_2 \oplus \u_0$, the structure matrix $A$ is of the form
\begin{align*}
    A = 
    \left[
	\begin{array}{c | c}      
		a & v^t \\
            \hline
            0 & A_0
	\end{array} \right] \quad \text{with }  v\in \R^{2n-2}  \text{ and } 
    A_0 \in \mathfrak{sp}(\omega_0) \cong \mathfrak{sp}(n-1, \R).  
\end{align*}    
\end{enumerate}    
\end{proposition} 

\noindent Here we are using the following convention: if $(V,\omega)$ is a symplectic vector space then
\[ \mathfrak{sp}(\omega):=\{T\in \operatorname{End}(V) \mid \omega(Tx,y)+\omega(x,Ty)=0 \text{ for all } x,y\in V \}.\]
If $V=\R^{2d}$ and $\omega_{st}$ is the standard symplectic form on $\R^{2d}$, then we denote
\[ \mathfrak{sp}(d,\R):=\mathfrak{sp}(\omega_{st})\subset \mathfrak{gl}(2d,\R).\]
If the symplectic vector space $(V,\omega)$ has dimension $2d$ then, after fixing a symplectic basis of $V$, we obtain an isomorphism
\[ \mathfrak{sp}(\omega)\cong \mathfrak{sp}(d,\R). \]

\begin{remark}
    Following the description in Proposition \ref{prop: LW}(iii), $\g_A$ is unimodular if and only if $a=0$, since it is well known that matrices in $\mathfrak{sp}(d,\R)$ are traceless. 
\end{remark} 

\indent We collect some elementary results on exterior derivatives of forms on a symplectic almost abelian Lie algebra for future reference. For it we need to choose a basis $\{e_l \mid 1 \leq l \leq 2n-2\}$ of $\u_0$. 
\begin{lemma} \label{lemma: elementary facts}
    The following statements hold in a symplectic almost abelian Lie algebra $\g_A = \R f_1 \ltimes \u$ equipped with a symplectic form $\omega = f^1 \wedge f^2 + \omega_0$ as described in Proposition \ref{prop: LW}. 
\begin{enumerate} [\rm (i)]
    \item If $\eta \in \alt^k \g_A^\ast$ is exact then $\eta = f^1 \wedge \mu$ for some $\mu \in \alt^{k-1} \u^*$.
    \item If $\eta \in \alt^k \g_A^\ast$ is of the form $\eta = f^1 \wedge \mu$ for some $\mu \in \alt^{k-1} \g_A^*$ then it is closed. 
    \item If $\alpha \in \alt^2 \g_A^\ast$ is a closed $2$-form then $\alpha([f_1,x],y)=-\alpha(x,[f_1,y])$ for all $x,y \in \u$. 
    \item $df^1 = 0$ and $d f^2 = -f^1 \wedge (a f^2 + \nu)$ for $\nu = \sum_{i=1}^{2n-2} v_i e^i \in \u_0^*$. In particular, if $a = 0$ and $v = 0$ then $df ^2 = 0$. 
    \item If $a = 0$ then $f^1 \wedge f^2$ is closed and non-exact.      
\end{enumerate}    
\end{lemma}     
\begin{proof} 
    (i) is immediate from the from the bracket relations in equation \eqref{eq: bracket relations}, and (ii) follows directly from (i). Also, (iii) follows easily from equation \eqref{eq: d2}. Certainly $f ^1$ is closed, and the rest of (iv) is established in a straightforward manner from Proposition \ref{prop: LW}. Lastly, for (v), $f^1 \wedge f^2$ is closed from (ii); and if it were exact, say $f^1 \wedge f^2 = d \eta$ for some $\eta = c_1 f^1 + c_2 f^2 + \theta$ with $\theta \in \u_0^*$, then (i) and (iv) entail the absurd that 
\begin{align*}
    f^1 \wedge f^2 &= c_2 df^2 + d \theta = f^1 \wedge \mu \text{ for some $\mu \in \u_0^*$}. \qedhere 
\end{align*}  
\end{proof}
\indent The following general results are used in Section \ref{section: the space of symplectic forms}.  
\begin{lemma} \label{lemma: equivalencia para morfismos}
    A linear map $f:\g \to \g$ is a Lie algebra homomorphism if and only if its pullback $f^*: \alt^* \g^* \to \alt^* \g^*$ satisfies $d(f^*\alpha)=f^*(d\alpha)$ for all $\alpha \in \g^*$. 
\end{lemma}
\begin{proof}
    It follows directly from the validity of the next two equations:
\begin{gather*}
    d(f^* \alpha)(x,y) = - (f^*\alpha)([x,y]) = - \alpha( f [x,y] ), \\
    f^*(d \alpha)(x,y) = (d \alpha)(f(x), f(y)) = -\alpha( [f(x), f(y)]),
\end{gather*} 
    \noindent for any $x$, $y \in \g$.
\end{proof}
\begin{lemma} \label{lemma: propagacion algebraica}
    Let $\g$ be a Lie algebra equipped with a symplectic form $\omega$. Assume that $\g$ has a decomposition as a semidirect product $\g=\h_1 \ltimes \h_2$ such that $\h_1$ and $\h_2$ are $\omega$-orthogonal subspaces. Then both the subalgebra $\h_1$ and the ideal $\h_2$ are non-degenerate with respect to $\omega$, and in particular $\omega_1:=\omega|_{\h_1\times \h_1}$ is a symplectic form on $\h_1$. Moreover, if the Lefschetz condition fails in $(\h_1,\omega_1)$ at degree $k$ then it also fails for $(\g,\omega)$ at degree $k$.
\end{lemma}
\begin{proof}
    The first statement is clear. For the second, note that the projection $\pi:\g\to \h_1$ is a surjective Lie algebra homomorphism that induces an injective linear map in cohomology, $\iota: H^k(\h_1)\to H^k(\g)$. The proof then follows the same lines as in the proof of Lemma \ref{lemma: propagacion}.
\end{proof}
\begin{remark}
    It is tempting to conclude that the key point in Lemma \ref{lemma: propagacion algebraica} lies not in the splitting of the symplectic Lie algebra $(\g, \omega)$ as a semidirect product in a way compatible with $\omega$, but rather in the splitting of a symplectic form such that one of its factors is, in a sense, non-hard-Lefschetz. However, this intuition is incorrect: without additional hypotheses ensuring that the inclusion $H^k(\h_1) \to H^k(\g)$ is injective, there is no reason to expect that the Lefschetz operator at degree $k$ has a nontrivial kernel. An explicit counterexample appears in \cite[Section 2]{Sawai}, and it is of the form $\g = \R^2 \ltimes \n$, where $\n$ is a six-dimensional nilpotent subalgebra of $\g$. Therein, it is shown that $(\g, \omega)$ satisfies the hard-Lefschetz condition for some symplectic form $\omega$ with respect to which $\R^2$ and $\n$ are $\omega$-orthogonal. This is at odds with the fact that $(\n, \omega \vert_{\n \times \n})$ does not satisfy the hard-Lefschetz condition because of Theorem \ref{thm: Benson y Gordon}, since $\n$ is nilpotent.
\end{remark} 
\indent The class of manifolds considered in this article exhibit the behavior described in Lemma \ref{lemma: propagacion algebraica}, and it is this feature that ultimately leads to our main result.

\indent It is not surprising that the Jordan canonical form of the structure matrix $A$ of the Lie algebras under consideration plays a role in characterizing which ones satisfy the hard-Lefschetz condition. With this in mind, we recall the following result.
\begin{proposition} \cite{etal} \label{prop: laura 1}
    Fix a scalar $b \in \R$, vectors $w \in \R^k$ and $u \in \R^m$, and matrices $Q \in M_k(\R)$ and $R \in M_m(\R)$. Suppose additionally that $b \notin \spec(Q)$.
\begin{enumerate} [\rm (i)]
    \item The following matrices are conjugate: 
\begin{align*} 
    \left[
	\begin{array}{c | c}      
		b & w^t \\
            \hline 
            0 & Q
	\end{array} \right] \quad \text{and} \quad 
    \left[
	\begin{array}{c | c}      
		b & 0 \\
            \hline
            0 & Q
	\end{array} \right].
\end{align*}
    \item The following matrices are conjugate:
\begin{align*} 
    \left[
	\begin{array}{c | c | c}      
		b & u^t & w^t \\
            \hline 
            0 & R & 0 \\
            \hline 
            0 & 0 & Q
	\end{array} \right] \quad \text{and} \quad 
    \left[
	\begin{array}{c | c | c}      
		b & u^t & 0\\
            \hline 
            0 & R & 0 \\
            \hline
            0 & 0 & Q
	\end{array} \right].
\end{align*}
\end{enumerate}
\end{proposition}


\smallskip

\indent For those matrices under our scrutiny,
\begin{align*}
    A = 
    \left[
	\begin{array}{c | c}      
		a & v^t \\
            \hline
            0 & A_0
	\end{array} \right] \quad \text{with} \quad
    A_0 \in \mathfrak{sp}(n-1, \R), 
\end{align*}
\noindent Proposition \ref{prop: laura 1} shows that the Jordan canonical form of $A$ is partially determined by that of $A_0$, the only complication arising when $a \in \spec(A_0)$. Proposition \ref{prop: laura 2} below deals with the possible Jordan canonical forms of matrices belonging to $\mathfrak{sp}(m,\R)$ that possess only real eigenvalues (see also Remark \ref{obs: laura 2 es mas general}). We need to introduce some notation in order to state the result precisely. 

\indent Denote by $J_m(\lambda)$ the elementary Jordan block of size $m \times m$ with eigenvalue $\lambda \in \mathbb{R}$,  
\begin{align} \label{eq: elementary jordan block}
    J_m(\lambda) = 
    \left[
	\begin{array}{c c c c c c}      
		\lambda & 0 & 0 & \cdots & 0 & 0\\ 
        1 & \lambda & 0 & \cdots & 0 & 0\\ 
        0 & 1 & \lambda & \cdots & 0 & 0 \\
        \vdots & \vdots & \vdots & & \vdots & \vdots \\
        0 & 0 & 0 & \cdots & \lambda & 0 \\
        0 & 0 & 0 & \cdots & 1 & \lambda
	\end{array} \right].
\end{align}
\noindent We know by elementary linear algebra that all matrices $M \in M_m(\R)$ with a unique real eigenvalue $\lambda$ are either $\lambda \I_{m \times m}$ or conjugate to
\begin{align*}
    J_{n_1}(\lambda)^{p_1} \oplus \cdots \oplus J_{n_k}(\lambda)^{p_k} \oplus\lambda \I_{t \times t}, \quad k \geq 1, \; n_1 > \cdots > n_k \geq 2, \; p_i\geq 1 \text{ for all } i, \; t \geq 0.  
\end{align*}
\noindent In particular, every such matrix $M$ has an associated tuple
\begin{align*}
    \Lambda(M) = (n_1, \ldots, n_k; p_1, \ldots, p_k; t), \quad \text{where} \quad p_1 n_1 + \cdots + p_k n_k + t = m,
\end{align*}
\noindent representing its Jordan canonical form. Moreover, if $M = M_1 \oplus \dots \oplus M_s$ with $\spec(M_j) = \{\lambda_j\}$ and $\lambda_j \neq \lambda_i$ for all $i \neq j$ then we can associate multiple tuples to $M$, namely
\begin{align*}
    \Lambda(M_1), \quad \ldots, \quad \Lambda(M_s). 
\end{align*} 

\noindent The following result was originally proved in \cite{MMRR} but we use the formulation given in \cite{etal}, which is more suitable to our purposes. 

\begin{proposition} \cite{etal} \label{prop: laura 2}
    Let $M \in M_{2m}(\R)$ with distinct real eigenvalues $\{\lambda_1, \ldots, \lambda_s\}$. Decompose $M = M_1 \oplus \cdots \oplus M_s$, where $\spec(M_j) = \{\lambda_j\}$ for all $1 \leq j \leq s$. Then $M$ is conjugate to a matrix in $\mathfrak{sp}(m, \R)$ if and only if the following conditions are both satisfied: 
\begin{enumerate} 
    \item if $\lambda_j = 0$ and $\Lambda(M_j) = (n_1, \ldots, n_k ; p_1, \ldots, p_k;t)$ then $p_r$ is even for all odd $n_r$ and $t$ is even,
    \item if $\lambda_j \neq 0$ then there exists $i \neq j$ such that $\lambda_i = - \lambda_j$ and $\Lambda(M_i) = \Lambda(M_j)$.
\end{enumerate}
\end{proposition}

\begin{remark} \label{obs: laura 2 es mas general}
    What is proven in \cite{etal} is more general than what was stated in Proposition \ref{prop: laura 2}, if anything because the case with complex eigenvalues is covered there. We do not consider this case since we are interested in Lie algebras $\g_A$ that are completely solvable, for which we need $A$ to have only real eigenvalues.
\end{remark}

\indent We are interested only in the completely solvable and unimodular case, so we assume throughout that $A$ (as in Proposition \ref{prop: LW}) has real eigenvalues and $a = 0$. In this scenario, Proposition \ref{prop: laura 1} implies that $A$ is then conjugate to either
\begin{align} \label{eq: the alternative}
    \left[
	\begin{array}{c | c}      
		0 & 0 \\
            \hline 
            0 & A_0
	\end{array} 
    \right]
        \quad \text{or} \quad 
    \left[
	\begin{array}{c | c | c}      
		0 & v & 0\\
            \hline 
            0 & Z & 0 \\
            \hline
            0 & 0 & N
	\end{array}
    \right] = 
    \left[
	\begin{array}{c | c }       
            M & 0  \\
            \hline
             0 & N
	\end{array}
    \right],
\end{align}
\noindent depending on whether $0 \notin \spec(A_0)$ or $0 \in \spec(A_0)$, respectively. Note that for $v$ to be nonzero it must happen that $0 \in \spec(A_0)$; however, it is possible for $v$ to be zero and that $0 \in \spec(A_0)$. For the case $0 \in \spec(A_0)$, we identify $A_0$ with $Z \oplus N$, so that $Z$ is the restriction of $A_0$ to the generalized eigenspace $V_0$ of $0 \in \spec(A_0)$ and $N$ is the restriction of $A_0$ to the direct sum of the other generalized eigenspaces of $A_0$ (more on that in the next section). From here on, we work exclusively with matrices as in equation \eqref{eq: the alternative}, and moreover we assume $A_0$ to be in Jordan basis. 

\section{The hard-Lefschetz condition on almost abelian Lie algebras} \label{section: main section}

\subsection{Special notations and nomenclature.} \label{section: notation}

In this section we establish notation and nomenclature that will be referred to throughout the remainder of the article. We use the notation in Proposition \ref{prop: LW} for a fixed symplectic form $\omega$ on $\g_A$.

\indent Let $\Xi := \{ \lambda_l \mid 1 \leq l \leq q \}$  denote the set of positive eigenvalues of $A_0$. For every eigenvalue $\lambda \in \spec(A_0)$ of $A_0$, not necessarily nonzero, denote by $V_{\lambda}$ the generalized eigenspace of $A_0$ associated to $\lambda$. Taking into account the results of Proposition \ref{prop: laura 2}, for all $1 \leq l \leq q$ we set $W_{\lambda_l} := V_{\lambda_l} \oplus V_{- \lambda_l}$. When there is no risk of confusion, we denote $W_l := W_{\lambda_l}$. We refer to these as the \textit{double generalized eigenspaces of $A_0$}. We also define $W_0 := \R f_2 \oplus V_0$, so as to have the following vector space decompositions:
\begin{align*}
    \u_0 = V_0 \oplus W_1 \oplus \cdots \oplus W_q, \quad \u = W_0 \oplus W_1 \oplus \cdots \oplus W_q.
\end{align*}
\noindent Certainly, as $\g_A = \R f_1 \oplus \u$,
\begin{align} \label{eq: g_A es f1 y u}
    \alt^{\ell} \g_A^* = \R f^1 \otimes \alt^{\ell - 1} \u^* \oplus \alt^{\ell} \u^*
\end{align}
\noindent and, for all $0 \leq \ell \leq \dim \g_A$,  
\begin{align} \label{eq: double generalized eigenspaces}
    \alt^{\ell} \u^* = \bigoplus_{i_0 + i_1 + \cdots + i_q = \ell} \alt^{i_0, i_1, \ldots, i_q}, \quad \text{where } \alt^{i_0, \ldots, i_q} := \alt^{i_0} W_0^* \otimes \alt^{i_1} W_1^* \otimes \cdots \otimes \alt^{i_t} W_q^*
\end{align} 

\indent Set $W := \bigoplus_{l=1}^q W_l$ to be the complementary subspace of $W_0$ in $\u$, so $\u = W_0 \oplus W$. The discussion leading up to equation \eqref{eq: the alternative} implies that a matrix $A$ with 
 $a = 0$ and possibly with $v \neq 0$ is then conjugate to either
\begin{align*}
    \left[
	\begin{array}{c | c}      
		0 & 0 \\
            \hline 
            0 & A_0
	\end{array} 
    \right]
        \quad \text{or} \quad 
    \left[
	\begin{array}{c | c | c}      
		0 & u & 0\\
            \hline 
            0 & Z & 0 \\
            \hline
            0 & 0 & N
	\end{array}
    \right] = 
    \left[
	\begin{array}{c | c }       
            M & 0  \\
            \hline
             0 & N
	\end{array}
    \right],
\end{align*}
\noindent depending on whether $0 \notin \spec(A_0)$ or $0 \in \spec(A_0)$, respectively. It is clear that
\begin{align*}
    M := A \vert_{W_0}, \quad Z := A_0 \vert_{W_0}, \quad N := A_0 \vert_W;  
\end{align*}
\noindent moreover, $\spec(Z) = \spec(M) = \{0\}$ and $0 \notin \spec(N)$. 

\indent The rest of the definitions in this section refer only to $A_0$ and $\u_0$. Pick an arbitrary eigenvalue $\lambda \in \spec(A_0)$ of $A_0$, possibly zero. As observed before, the restriction $A_0 \vert_{V_{\lambda}}$ is either $\lambda \I_{m \times m}$ or conjugate to 
\begin{align*}
    J_{n_1}(\lambda)^{p_1} \oplus \cdots \oplus J_{n_k}(\lambda)^{p_k} \oplus\lambda \I_{t \times t}, \quad k \geq 1, \; n_1 > \cdots > n_k \geq 2, \; p_i \geq 1 \text{ for all } i, \; t \geq 0.  
\end{align*}
\noindent In the latter case, there are subspaces $X_{n_1}(\lambda), \ldots, X_{n_k}(\lambda)$ of $V_{\lambda}$ such that 
\begin{align*}
    A_0 \vert_{X_{n_l} (\lambda)} \text{ and } J_{n_l}(\lambda) \text{ are conjugate for all $1 \leq l \leq k$};
\end{align*}
\noindent what we mean is that these subspaces correspond to elementary Jordan blocks of $A_0$, and thus we call them \textit{elementary canonical subspaces} of $A_0$. We also say that, for all $1 \leq l \leq k$, $X_{n_l}$ has \textit{size $n_l$}. We set
\begin{align*}
    X_{n_l}(\lambda) = \Span_{\mathbb{R}} \{x_i \mid 1 \leq i \leq n_l\},
\end{align*}
\noindent and consistently think of $\{x_i \mid 1 \leq i \leq n_l\}$ as a Jordan basis, so $A_0 \vert_{X_{n_l} (\lambda)}$ and $J_{n_k}(\lambda)$ are actually equal for all $1 \leq l \leq k$. Referring again to Proposition \ref{prop: laura 2}, if $\lambda \in \spec(A_0)$ is nonzero then $- \lambda \in \spec(A_0)$ is also a nonzero eigenvalue of $A_0$, and moreover for every elementary canonical subspace $X_r(\lambda) \subseteq V_{\lambda}$ of size $r$ of $A_0$ there exists an elementary canonical subspace $X_r(- \lambda) \subseteq V_{- \lambda}$ of size $r$ of $A_0$. Making non-canonical identifications, one can speak of ``the"\! companion subspace $X_r(- \lambda)$ of $X_r(\lambda)$. In this regard, we denote the basis of each of them as
\begin{align} \label{eq: x tilde base} 
    X_r(\lambda) = \Span_{\mathbb{R}} \{x_i \mid 1 \leq i \leq n_l\}, \quad X_r(- \lambda) = \Span_{\mathbb{R}} \{x_{i+r} \mid 1 \leq i \leq n_l\},
\end{align}
\noindent When there is no risk of confusion, we call $X_r( \lambda)$ of $X_r(-\lambda)$ respectively by $X^+$ and $X^-$. We say that $\widetilde{X} := X^+ \oplus X^-$ is a  \textit{double elementary canonical subspace of size $r$}  associated with a pair of nonzero eigenvalues $\{\lambda, - \lambda \} \subseteq \spec(A_0)$ of $A_0$. Note that we say that $\widetilde{X}$ has \textit{size} $r$ even though it has \textit{dimension} $2r$. Associated to every double eigenspace $W_i$ ($i>0$) there exists a sequence of double elementary canonical subspaces $\widetilde{X_{i,1}}, \ldots, \widetilde{X_{i,p_i}}$, possibly of repeating sizes, such that 
\begin{align} \label{eq: a finer decomposition}
    W_i = \widetilde{X_{i,1}} \oplus \cdots \oplus \widetilde{X_{i,p_i}}.
\end{align}

\subsection{The space of closed 2-forms.} \label{section: the space of closed 2 forms}

Consider a matrix $A$ as in equation \eqref{eq: the alternative}, with $a = 0$ and possibly $v \neq 0$. Refer to the previous Section for the notation. Our first intent is to describe the behavior in cohomology of the symplectic form $\omega$ on $\g_A$. As it turns out, for the case $v = 0$ (in equation \eqref{eq: the alternative}), we are able to provide a comprehensive description of the second-degree cohomology $H^2(\g_A)$ of $\g_A$. 

\indent As our first step, for a nonzero eigenvalue $\lambda \in \spec(A_0)$ of $A_0$, pick an elementary canonical subspace $X_r(\lambda)$ of size $r$. As there is no risk of confusion, set $X := X_r(\lambda)$ from this section on. Recall, 
\begin{align*}
    X = \Span_{\mathbb{R}} \{x_j \mid 1 \leq j \leq r \} \subseteq V_{\lambda},
\end{align*}
\noindent and we assume that the basis for $X$ is the Jordan basis. This last bit and the fact that $\lambda \neq 0$ amount to the validity of the following relations:
\begin{align} \label{eq: bracket relations X}
    [f_1, x_j] = 
    \begin{cases}
        \lambda x_j + x_{j+1}, & j < r, \\
        \lambda x_r, & j = r. 
    \end{cases} 
\end{align}
\noindent Note that, if $r = 1$, these equations reduce to $[f_1,x_1] = \lambda x_1$. Proceeding recursively from here, one finds that
\begin{align} \label{eq: vj} 
    x_j = \left[f_1, \frac{1}{\lambda} x_j + \sum_{i = j+1}^r \frac{(-1)^{i-j}}{\lambda^{i-j+1}} x_i \right] =: \left[f_1, \frac{1}{\lambda} x_j + \sum_{i = j+1}^r c_{i,j} x_i \right] \quad \text{ for all $1 \leq j \leq r$}.
\end{align}
\noindent We have set $c_{i,j} := \frac{(-1)^{i-j}}{\lambda^{i-j+1}}$ for simplicity. It is important to keep in mind the following relations for the $c_{i,j}$'s, valid for all indices that make sense: 
\begin{align} \label{eq: cij}
    c_{i,j} \neq 0, \quad \lambda^2 c_{j+1,j} = -1, \quad \lambda c_{i+1,j} + c_{i,j} = 0. 
\end{align}
\noindent Now pick another eigenvalue $\mu \in \spec(A_0)$, possibly zero. Also, pick an elementary canonical subspace $Y := Y_s(\mu)$ of size $s$, and again assume that the basis $\{y_i \mid 1 \leq i \leq s\}$ for $Y$ is the Jordan basis. We are \textit{not} assuming that $r = s$. If $\mu \neq 0$ then similar relations to those in equation \eqref{eq: bracket relations X} hold; in the special case that $\mu = 0$, the following relations hold: 
\begin{align*}
    [f_1, y_k] = 
    \begin{cases}
        v_k f_2 + y_{k+1}, & k < s, \\
        v_s f_2, & k = s, 
    \end{cases} 
\end{align*}
\noindent where $\sum_{k=1}^s v_k y_k$ is the projection of $v$ in equation \eqref{eq: the alternative} onto $Y$. 

\indent It turns out that much can be said about the restriction to either $X$, $Y$, or $X \oplus Y$ to either a generic closed $2$-form $\alpha$ on $\g_A$ when $v = 0$, or for the symplectic form $\omega$ for $v \neq 0$.
\begin{lemma} \label{lemma: X, Y uno}
    Fix $\lambda$, $\mu \in \spec(A_0)$ as above, with $\lambda \neq 0$ and $\mu \neq - \lambda$. Let $\omega \in \alt^2 \g_A^*$ be a symplectic form as in Proposition \ref{prop: LW}. 
\begin{enumerate} [\rm (i)]
    \item If $\mu \neq 0$ then $X$ and $Y$ are $\alpha$-orthogonal to each other for any closed $2$-form $\alpha \in \alt^2 \g_A^{\ast}$.
    \item If $\mu = 0$ and $v = 0$ then $X$ and $Y$ are $\alpha$-orthogonal to each other for any closed $2$-form $\alpha \in \alt^2 \g_A^{\ast}$.
    \item If $\mu = 0$ and $v \neq 0$ then $X$ and $Y$ are $\omega$-orthogonal.
    \item $X$ is $\alpha$-isotropic for any closed $2$-form $\alpha \in \alt^2 \g_A^{\ast}$.  
\end{enumerate}
\end{lemma}
\begin{proof} 
    We begin by establishing (i) and (ii), which can be done simultaneously. For (iii), it can be proven along similar lines, with the additional caveat that $f_2$ is $\omega$-orthogonal to any element in $\u_0$ as per Proposition \ref{prop: LW}. The validity of (iv) results directly from (i), taking $\mu = \lambda$ (thus $\mu\neq -\lambda$) and $Y = X$, since then we obtain that $X$ is $\alpha$-orthogonal to itself.
    
    \indent So let's prove (i) and (ii). Notice that $\alpha(x_r, y_s) = 0$. Indeed, using Lemma \ref{lemma: elementary facts}(iii) and equation \eqref{eq: vj}, we get
\begin{align*}
    \alpha(x_r, y_s) = \frac{1}{\lambda} \alpha([f_1, x_r], y_s) = -\frac{1}{\lambda} \alpha(x_r, [f_1, y_s]) = - \frac{\mu}{\lambda}\alpha(x_r, y_s),
\end{align*}
    \noindent and so the relation $\alpha(x_r, y_s) = 0$ follows from $\mu \neq -\lambda$. Thus, there exists a minimal $j_0 \leq r$ such that $\alpha(x_j, y_s) = 0$ for all $j_0 \leq j \leq r$. It must be the case that $j_0 = 1$, for otherwise $j_0 > 1$ and, from Lemma \ref{lemma: elementary facts}(iii) and equation \eqref{eq: vj} again, we obtain
\begin{align*} 
    \alpha(x_{j_0-1}, y_s) &= \alpha \left( \left[f_1, \frac{1}{\lambda} x_{j_0-1} + \sum_{i = j_0}^r c_{i,j_0 -1} x_i \right], y_s \right) = - \alpha \left( \frac{1}{\lambda} x_{j_0-1} + \sum_{i = j_0}^r c_{i, j_0 -1} x_i, [f_1, y_s] \right) \\
    &= -\mu \alpha\left(\frac{1}{\lambda}x_{j_0-1} + \sum_{i = j_0}^r c_{i, j_0 -1} x_i, y_s \right) = - \frac{\mu}{\lambda}\alpha(x_{j_0-1}, y_s),
\end{align*}
    \noindent which implies that $\alpha(x_{j_0-1}, y_s) = 0$ since $\mu\neq -\lambda$, contradicting the minimality of $j_0$. Thus, there exists a minimal $k_0 \leq s$ such that $\alpha(x_j, y_k) = 0$ for all $1 \leq j \leq r$ and all $k_0 \leq k \leq s$. Again, it must be the case that $k_0 = 1$, for otherwise $k_0 > 1$ and, from Lemma \ref{lemma: elementary facts}(iii) and equation \eqref{eq: vj} once more, we obtain
\begin{align*} 
    \alpha(x_j,y_{k_0-1}) & =  \alpha\left( \left[f_1,\frac{1}{\lambda}x_j + \sum_{i = j+1}^r c_{i,j} x_i\right],y_{k_0-1}\right)=-\alpha\left(\frac{1}{\lambda}x_{j}+ \sum_{i = j+1}^r c_{i,j} x_i,[f_1,y_{k_0-1}]\right) \\
    &= - \alpha\left(\frac{1}{\lambda}x_{j}+ \sum_{i = j+1}^r c_{i,j} x_i, \mu y_{k_0-1} + y_{k_0}) \right) = -\frac{\mu}{\lambda}\alpha(x_{j},y_{k_0-1}),
\end{align*}
    \noindent which implies that $\alpha(x_j, y_{k_0-1}) = 0$ for any $j$ since $\mu\neq -\lambda$, contradicting the minimality of $k_0$. Thus, $\alpha(x,y) = 0$ for all $x \in X$ and $y \in Y$, and therefore $X$ and $Y$ are $\alpha$-orthogonal to each other.  
\end{proof}
    \indent Note that Lemma \ref{lemma: X, Y uno}(i) also implies that $X$ and $Y$ are $\alpha$-orthogonal to each other for any closed $2$-form $\alpha$ on $\g_A$ even when both $X$ and $Y$ are associated to the same eigenvalue $\lambda$.   
    
    \indent For the next result, set $m := \min \{r,s\}$. Note that the statements are meaningless for the case $r = s = 1$, corresponding to the ``diagonalizable bits"\! within the Jordan decomposition of $A_0$: nothing can be said about the matrix form of $\alpha$ in that case. Similarly, when either $r = 1$ or $s = 1$ then only the first statement is meaningful. As this result does not involve $V_0$, there is no need to distinguish between the $v = 0$ and $v \neq 0$ cases. 
\begin{lemma} \label{lemma: X, Y dos}
    In the notation described above, assuming additionally that $\mu = - \lambda$ and both nonzero, the following statements hold for any closed $2$-form $\alpha \in \alt^2 \g_A^{\ast}$: 
\begin{enumerate} [\rm (i)]
    \item $\alpha(x_j, y_k) = 0$ for all $m +2 \leq j+k \leq r + s$.
    \item $\alpha(x_{j+1}, y_k) + \alpha(x_j, y_{k+1}) = 0$ for all $1 \leq j \leq r - 1$ and $1 \leq k \leq s - 1$.
\end{enumerate} 
\end{lemma}
\begin{proof} 
    \noindent Both statements are proven at the same time. By interchanging the roles of $X$ and $Y$ and renaming if necessary, there is no loss in generality in assuming that $s \leq r$, and thus $m = s$. We follow a similar argument as in the proof of Lemma \ref{lemma: X, Y uno}.

    \indent For all $1 \leq j \leq r - 1$, by Lemma \ref{lemma: elementary facts}(iii) and equation \eqref{eq: vj}, we have that
\begin{align*}
    \alpha(x_j, y_s) &= \alpha \left( \left[ f_1, \frac{1}{\lambda} x_j + \sum_{i = j+1}^r c_{i,j} x_i \right], y_s \right) = - \alpha \left( \frac{1}{\lambda} x_j + \sum_{i = j+1}^r c_{i,j} x_i , \left[ f_1, y_s \right] \right) \\
    &= - \alpha \left( \frac{1}{\lambda} x_j + \sum_{i = j+1}^r c_{i,j} x_i , - \lambda y_s \right) = \alpha(x_j, y_s) + \lambda \alpha \left(\sum_{i = j+1}^r c_{i,j} x_i, y_s \right), 
\end{align*}
    \noindent and thus, as $\lambda \neq 0$,
\begin{align*}
    0 = \sum_{i = j+1}^r c_{i,j} \alpha \left( x_i, y_s \right) \text{ for all $1 \leq j \leq r-1$}.
\end{align*}    
    \noindent This is a upper-triangular linear system of equations for the unknowns $\alpha(x_2, y_s), \ldots, \alpha(x_r, y_s)$, and therefore can be solved iteratively, in a bottom-up fashion, substituting each found unknown in the next equation up. The first equation, indexed by $j = r-1$, is simply $0 = c_{r, r-1} \alpha(x_r, y_s)$. It follows that $\alpha(x_r, y_s) = 0$ since $c_{r,r-1}$ is nonzero (see equation \eqref{eq: cij}), and from there that 
\begin{align*}
    \alpha(x_r, y_s) = \cdots = \alpha(x_2, y_s) = 0. 
\end{align*}
    \noindent No information is found for $\alpha(x_1, y_s)$. The rest of the proof goes on from here following a routine inductive procedure that we describe as follows. Similarly as above, for all $1 \leq j \leq r - 1$ and $1 \leq k \leq s-1$, again by Lemma \ref{lemma: elementary facts}(iii) and equation \eqref{eq: vj}, we have that
\begin{align*}
    \alpha(x_j, y_k) &= \alpha \left( \left[ f_1, \frac{1}{\lambda} x_j + \sum_{i = j+1}^r c_{i,j} x_i \right], y_k \right) = - \alpha \left( \frac{1}{\lambda} x_j + \sum_{i = j+1}^r c_{i,j} x_i , \left[ f_1, y_k \right] \right) \\
    &= - \alpha \left( \frac{1}{\lambda} x_j + \sum_{i = j+1}^r c_{i,j} x_i, - \lambda y_k + y_{k+1} \right) \\
    &= \alpha(x_j, y_k) - \frac{1}{\lambda} \alpha \left( x_j, y_{k+1} \right) + \lambda  \sum_{i = j+1}^r c_{i,j} \alpha(x_i, y_k) - \sum_{i = j+1}^r c_{i,j} \alpha(x_i, y_{k+1}) 
\end{align*}
    \noindent The idea is to apply backwards induction in $k$ but only in statement (i): at every step, statement (i) is established by considering all $j$ such that $j + k \geq s +2$, while statement (ii) is established by considering all the other $j$. Note that the case $k = s$ is already covered. Thus, assume that the inductive hypothesis holds for values greater than $k$, meaning that 
\begin{gather} \label{eq: HI 1}
    \alpha(x_j, y_{k+1}) = \cdots = \alpha(x_r, y_{k+1}) = 0 \text{ for all $j > s - k$}.
\end{gather}    
    \indent Let's focus on the case $j > s-k$ first. After canceling out the $\alpha(x_j, y_k)$'s in the foregoing equation and using the first part of equation \eqref{eq: HI 1}, we arrive at the following:
\begin{align*} 
    0 = \sum_{i = j+1}^r c_{i,j} \alpha(x_i, y_k)  \text{ for all $j > s-k$}.
\end{align*}
    \noindent As before, for $j > s-k$, we have a upper-triangular linear system of equations for the unknowns $\alpha(x_{s-k+2}, y_k), \ldots, \alpha(x_r, y_k)$, from where it is direct to conclude that
\begin{align} \label{eq: HI 2}
    \alpha(x_r, y_k) = \cdots = \alpha(x_{s-k+2}, y_k) = 0.
\end{align}
    \noindent This fact alone establishes statement (i). Now let's rework the system we are working with, again canceling out the  $\alpha(x_j, y_k)$'s and using in particular equations \eqref{eq: cij}, \eqref{eq: HI 1}, and \eqref{eq: HI 2}:
\begin{align*} 
    0 &= \alpha \left( x_j, y_{k+1} \right) - \lambda^2  \sum_{i = j+1}^r c_{i,j} \alpha(x_i, y_k) + \lambda \sum_{i = j+1}^r c_{i,j} \alpha(x_i, y_{k+1}) \\
    &= \alpha \left( x_j, y_{k+1} \right) - \lambda^2 c_{j+1, j} \alpha(x_{j+1}, y_k) - \lambda \sum_{m = j+1}^{s-k} \left[ \lambda c_{m+1,j} \alpha( x_{m+1}, y_k) - c_{m,j} \alpha(x_m, y_{k+1})  \right] \\
    &\phantom{\alpha \left( x_j, y_{k+1} \right) -} \; \, - \lambda^2 \sum_{i = s - k+2 }^r c_{i,j} \alpha(x_i, y_k) + \lambda \sum_{i = s - k + 1}^r c_{i,j} \alpha(x_i, y_{k+1}) \\ 
    &= \alpha \left( x_j, y_{k+1} \right) + \alpha(x_{j+1}, y_k) + \lambda \sum_{m = j+1}^{s-k} c_{m,j} [ \alpha(x_m, y_{k+1}) + \alpha( x_{m+1}, y_k)] 
\end{align*}
    \noindent This last system of equations is valid for all $j \leq s - k$, and once again it can be solved iteratively. We describe briefly how: When $j = s - k$ the sum is empty, so the equation we get is $\alpha(x_{s-k}, y_{k+1}) + \alpha( x_{s-k+1}, y_k) = 0$. This is one of the equations we are seeking to establish. When $j = s-k-1$, using equation \eqref{eq: cij} again, we arrive at   
\begin{align*}
     0 = \alpha \left( x_{s-k-1}, y_{k+1} \right) + \alpha(x_{s-k+1}, y_k) - \underbrace{[ \alpha(x_{s-k}, y_{k+1}) + \alpha( x_{s-k+1}, y_k) ]}_{\bullet};
\end{align*}
    \noindent the equation for $j = s-k$ establishes that the term $\bullet$ is identically zero, so we get \[ \alpha \left( x_{s-k-1}, y_{k+1} \right) + \alpha(x_{s-k+1}, y_k) = 0,\] which is another equation we are seeking to establish. We see thus that every term in the sum is zero as it corresponds to an equation established in a previous step. This concludes the proof. 
\end{proof} 
    
    \indent Note that Lemma \ref{lemma: X, Y dos} also implies that, for all closed $2$-forms $\alpha$, the scalars $(-1)^j \alpha(x_j, y_k)$, where $j+k$ is constant, are all equal. This roughly means that in matrix form, the coefficients of $\alpha$ share the same absolute value and alternate in sign along \textit{antidiagonals} of said matrix. The first part of the same Lemma ensures that actually many of said antidiagonals are identically zero. 
    
    \indent Let's shed some light on the Lemmas \ref{lemma: X, Y uno} and \ref{lemma: X, Y dos}. First, applying Lemma \ref{lemma: X, Y uno} to every elementary Jordan block within $V_{\lambda}$ and $V_{\mu}$, we obtain the following result. 
\begin{corollary} \label{cor: splitting}
    Fix $\lambda$, $\mu \in \spec(A_0)$ as above, with $\lambda \neq 0$ and $\mu \neq - \lambda$. Let $\omega \in \alt^2 \g_A^*$ be a symplectic form as in Proposition \ref{prop: LW}. 
\begin{enumerate} [\rm (i)]
    \item If $\mu \neq 0$ then $V_{\lambda}$ and $V_{\mu}$ are $\alpha$-orthogonal for all closed $2$-forms $\alpha \in \alt^2 \g_A^*$.
    \item If $\mu = 0$ and $v = 0$ then $V_{\lambda}$ and $V_{\mu}$ are $\alpha$-orthogonal for all closed $2$-forms $\alpha \in \alt^2 \g_A^*$.
    \item If $\mu = 0$ and $v \neq 0$ then $V_{\lambda}$ and $V_{\mu}$ are $\omega$-orthogonal. 
    \item $V_{\lambda}$ is $\alpha$-isotropic for all nonzero $\lambda \in \spec(A_0)$.  
\end{enumerate}     
\end{corollary}     
    \indent Turning now to Lemma \ref{lemma: X, Y dos}. We illustrate the facts established therein with two examples.
\begin{example}
    Suppose $A_0$ in equation \eqref{eq: the alternative} is conjugate to a matrix of the form
\begin{align*}
    J_4(\lambda) \oplus J_4(- \lambda),
\end{align*}
    \noindent where $J_4$ is as in equation \eqref{eq: elementary jordan block} with $\lambda\neq 0$. According to Corollary \ref{cor: splitting}, the matrix corresponding to any closed $2$-form $\alpha$ on $\g_A$, when restricted to the subspace spanning $A_0$ and written in the Jordan basis, is of the form
\begin{align*}
    \left[
	\begin{array}{c | c}      
		\textbf{0}_{4 \times 4} & - P \\
            \hline 
            P & \textbf{0}_{4 \times 4}
	\end{array} \right],
\end{align*}
    \noindent for some $P \in M_4(\R)$. Moreover, according to Lemma \ref{lemma: X, Y dos}, $P$ is an ``upper-antidiagonal alternating"\! Hankel matrix. This means 
    that it has one of the following forms: 
\begin{align*}
    \left[
	\begin{array}{c c c c}      
		a & - b & c & - d \\
        b & - c & d & \phantom{+} 0 \\
        c & - d & 0 & \phantom{+} 0 \\
        d & \phantom{+} 0 & 0 & \phantom{+} 0
	\end{array} \right] \quad \text{ or } \quad
    \left[
	\begin{array}{c c c c}      
		\phantom{+} a & \phantom{+} b & \phantom{+} c & \phantom{+} d \\
        - b & - c & - d & \phantom{+} 0 \\
        \phantom{+} c & \phantom{+} d & \phantom{+} 0 & \phantom{+} 0 \\
        - d & \phantom{+} 0 & \phantom{+} 0 & \phantom{+} 0
	\end{array} \right]. 
\end{align*} 
\end{example}
\begin{example}
    Suppose $A_0$ in equation \eqref{eq: the alternative} is conjugate to a matrix of the form
\begin{align*}
    J_4(\lambda) \oplus J_4(- \lambda) \oplus J_3(\lambda) \oplus J_3(- \lambda),
\end{align*}
    \noindent where $J_4$ and $J_3$ are as in equation \eqref{eq: elementary jordan block} with $\lambda\neq 0$. According to Corollary \ref{cor: splitting}, the matrix corresponding to any closed $2$-form $\alpha$ on $\g_A$, when restricted to the subspace spanning $A_0$ and written in the Jordan basis, is of the form
\begin{align*}
    \left[
	\begin{array}{c | c | c | c }      
            \textbf{0}_{4 \times 4} & - P & \textbf{0}_{4 \times 3} & -Q \\ 
            \hline 
            P & \textbf{0}_{4 \times 4} & -R & \textbf{0}_{3 \times 3} \\ 
            \hline
            \textbf{0}_{3 \times 4} & R & \textbf{0}_{3 \times 3} & -S \\
            \hline 
            Q & \textbf{0}_{3 \times 3} & S & \textbf{0}_{3 \times 3} \\
	\end{array} \right], 
\end{align*}
    \noindent for some matrices $P$, $Q$, $R$, and $S$ of respective sizes $4 \times 4$, $3 \times 3$, $4 \times 3$, and $3 \times 3$, one example of which given by
\begin{gather*}
    P = \left[
	\begin{array}{c c c c }      
		a & - b & c & - d  \\
            b & - c & d & \phantom{+} 0 \\
            c & - d & 0 & \phantom{+} 0 \\
            d & \phantom{+} 0 & 0 & 0
	\end{array} \right], \quad 
    Q = \left[
	\begin{array}{c c c c}      
		i & - j & k & \phantom{+} 0 \\
            j & - k & 0 & \phantom{+} 0 \\
            k & \phantom{+} 0 & 0 & \phantom{+} 0
	\end{array} \right], \\
    R = \left[ \begin{array}{c c c}    
            e & - f & g \\
            f & - g & 0 \\
            g & \phantom{+} 0 & 0 \\
            0 & \phantom{+} 0 & 0
        \end{array} \right], \quad
    S = \left[ \begin{array}{c c c}    
            l & - m & n \\
            m & - n & 0 \\
            n & \phantom{+} 0 & 0 
        \end{array} \right].
\end{gather*}
    \noindent Other possible admissible matrices $P$, $Q$, $R$, and $S$ must display a similar pattern of zeros and their coefficients should alternate signs in every nonzero antidiagonal. As in the previous examples, they are ``upper-antidiagonal alternating"\! Hankel matrices of sorts.
\end{example}
\indent Examples of larger dimension become unmanageable quickly. Hopefully these two are enough to get the gist of it. At any rate, they pave a way to describe the second-degree cohomology $H^2(\g_A)$ of $\g_A$.

\indent We close this section with a description of the non-exact $2$-forms arising from the wedge product between $f^1$ and some $\mu \in \g_A^*$. When $a = 0$ and $0 \notin \spec(A_0)$, we find that essentially the only one is $f^1 \wedge f^2$.

\begin{lemma} \label{lemma: gran alegria}
    In the notation described above, assuming additionally that $0 \notin \spec(A_0)$, if $\nu \in \u_0^*$ is a nonzero $1$-form then $f^1 \wedge \nu$ is exact. 
\end{lemma}
\begin{proof} 
    \noindent Fix a basis $\{e_l \mid 1 \leq l \leq m\}$ of $\u_0$, where $n = m + 1$. Assume first that $\u_0$ is a double elementary canonical subspace of $A_0$; that is, $\spec(A_0) = \{\lambda, - \lambda\}$ (with $\lambda > 0$) and $W = \widetilde{X}$. Assume further that the chosen basis is a nice basis for $\widetilde{X}$: what we mean is that its first half is the Jordan canonical basis associated to $X^+$, and that its second half is the Jordan canonical basis associated to $X^-$. Write $\nu = \sum_{j=1}^{2m} c_j e^j \in \u_0^*$, and define $\nu_1 := \sum_{j=1}^m c_j e^j$ and $\nu_2 := - \sum_{j=m+1}^{2m} c_j e^j$, so that $\nu = \nu_1 - \nu_2$. The fact that $\nu$ is nonzero entails that
\begin{align*}
    k_0 := \max \{ k \mid \text{ $c_k \neq 0$ or $c_{m+k} \neq 0$} \} \geq 1
\end{align*}
    \noindent is well defined. If $k_0 = 1$ then $\nu_1 = c_1 e^1$ and $\nu_2 = c_{m+1} e^{m+1}$, and thus 
\begin{gather*}
    d \nu_1 = - c_1 \lambda f^1 \wedge e^1, \quad d \nu_2 = - c_{m+1} \lambda f^1 \wedge e^{m+1},
\end{gather*}
    \noindent from where it follows that
\begin{align*}
    d (\nu_1 + \nu_2) = - \lambda f^1 \wedge (\nu_1 - \nu_2), 
\end{align*}
    \noindent thereby establishing that $f^1 \wedge \nu$ is exact. Now, if $k_0 \geq 2$,     
\begin{gather*}
    d \nu_1 = - f^1 \wedge \left( c_1 \lambda e^1 + \sum_{j=2}^{k_0} [c_j \lambda e^j + c_j e^{j-1}] \right), \\
    d \nu_2 = - f^1 \wedge \left( c_{m+1} \lambda e^{m+1} + \sum_{j=m+2}^{m + k_0} [ c_j \lambda e^j - c_j e^{j-1}] \right),
\end{gather*}
    \noindent so 
\begin{align*}
    d (\nu_1 + \nu_2) &= - \lambda f^1 \wedge (\nu_1 - \nu_2) - f^1 \wedge \left(\sum_{j = 2}^{k_0} c_j e^{j-1} - \sum_{j = m + 2}^{m + k_0} c_j e^{j-1} \right) \\
    &= - \lambda f^1 \wedge \nu - f^1 \wedge \left(\sum_{k = 1}^{k_0-1} c_{k+1} e^k - \sum_{k = m+1}^{m + k_0 - 1} c_{k+1} e^k \right). 
\end{align*}  
    \noindent Set $\mu_1 := \sum_{k = 1}^{k_0 - 1} c_{k+1} e^k$ and $\mu_2 := \sum_{k = m+1}^{m + k_0 - 1} c_{k+1} e^k$, so that 
\begin{align*}
    d (\nu_1 + \nu_2) = - \lambda f^1 \wedge \nu - f^1 \wedge (\mu_1 - \mu_2),
\end{align*}
    \noindent and $f^1 \wedge \nu$ is exact if and only if $f^1 \wedge (\mu_1 - \mu_2)$ is exact. After repeating this argument enough times, we see that $f^1 \wedge \nu$ is exact if and only if $f^1 \wedge ( c_{k_0} e^1 - c_{m+k_0} e^{m+1})$ is nonzero and exact. Since
\begin{align*}
    - \frac{1}{\lambda} d(c_{k_0} e^1 + c_{2 k_0} e^{m+1}) = f^1 \wedge \left( c_{k_0} e^1 - c_{2 k_0} e^{m+1} \right),
\end{align*}   
    \noindent we conclude that $f^1 \wedge ( c_{k_0} e^1 - c_{m+k_0} e^{m+1})$ is exact, and so is $f^1 \wedge \nu$.
    
    \indent For the general case, split $\nu$ as a sum $\theta_1 + \cdots + \theta_t$ arising from dualizing the decompositions
\begin{align*}
    \u_0 = V_0 \oplus W_1 \oplus \cdots \oplus W_q, \quad W_i = \widetilde{X_{i,1}} \oplus \cdots \oplus \widetilde{X_{i,p_i}} \text{ for all $1 \leq i \leq q$},
\end{align*}    
    \noindent and noting that $V_0=0$ because $0 \notin \spec(A_0)$. An argument along similar lines to the one described above works for every nonzero $\theta_j$.  
\end{proof} 
\begin{remark}
    The conclusion of Lemma \ref{lemma: gran alegria} is false for $0 \in \spec(A_0)$, even if $v = 0$. For example, for $A = (0) \oplus J_2(0) \oplus J_2(0)$ written in a basis $\{ f_2, e_1, e_2, e_3, e_4\}$, it is clear that $f^1 \wedge e^2$ is non-exact.
\end{remark}
\begin{corollary} \label{cor: para citar}
    If $0 \notin \spec(A_0)$ then, for $\mu \in \g_A^*$,  $f^1 \wedge \mu$ is non-exact if and only if $\mu = c_1 f^1 + c_2 f^2$ for some $c_1$, $c_2 \in \R$.    
\end{corollary}

\subsection{Computations in cohomology.} \label{section: computations in cohomology} 

Consider a matrix $A$ as in equation \eqref{eq: the alternative}, with $a = 0$ and possibly with $v \neq 0$. Refer to Section \ref{section: notation} for the notation; in particular, recall equation \eqref{eq: x tilde base} and the discussion leading up to it.  
\begin{definition} \label{def: circuits}
    Consider double elementary canonical subspaces $\widetilde{X}$ and $\widetilde{Y}$ of $A_0$, both associated to a pair of nonzero eigenvalues $\{\lambda, - \lambda\}$ of $A_0$, respectively of sizes $r$ and $s$. Set $m := \min\{r,s\}$. The \textit{circuits of $\widetilde{X}$ and $\widetilde{Y}$} are the following set of $2$-forms:
\begin{gather*} 
    g_l(x,x) := \sum_{i = 1}^l (-1)^{i + 1} x^i \wedge x^{ r + l + 1 - i }, \quad 1 \leq l \leq r, \\
    g_l(y,y) := \sum_{i = 1}^l (-1)^{i + 1} y^i \wedge y^{ s + l + 1 - i }, \quad 1 \leq l \leq s, \\  
    g_l(x,y) := \sum_{i = 1}^l (-1)^{i + 1} x^i \wedge y^{ m + l + 1 - i }, \quad 1 \leq l \leq m, \\ 
    g_l(y,x) := \sum_{i = 1}^l (-1)^{i + 1} y^i \wedge x^{ m + l + 1 - i }, \quad 1 \leq l \leq m.
\end{gather*}   
    \noindent The number of terms of a circuit is called its \textit{length}. 
\end{definition}
\indent This unorthodox notation seems to strike the best balance between precision and a not-so-cluttered appearance. Notice the ``$r$"\! in the definition of the $g_l(x,x)$'s, as well as the ``$s$"\! in the definition of the $g_l(y,y)$. These are replaced by ``$m$"\! in the definition of the $g_l(x,y)$'s and the $g_l(y,x)$'s. We list some circuits below to clarify what we mean:  
\begin{gather*}
    g_1(x,y) = x^1 \wedge y^{m+1}, \\
    g_2(x,y) = x^1 \wedge y^{m+2} - x^2 \wedge y^{m+1}, \\
    g_3(x,y) = x^1 \wedge y^{m+3} - x^2 \wedge y^{m+2} + x^3 \wedge y^{m+1}, \\
    g_4(x,y) = x^1 \wedge y^{m+4} - x^2 \wedge y^{m+3} + x^3 \wedge y^{m+2} - x^4 \wedge y^{m+1}, \\
    \vdots \\
    g_m(x,y) = x^1 \wedge y^{2m} - x^2 \wedge y^{2m-1} + \cdots + (-1)^{m+1} x^m \wedge y^{m+1}.
\end{gather*}
\indent If $\widetilde{X} = \widetilde{Y}$ then there are not as many circuits, since there cannot be mixing. In this case, all circuits are of the form $g_l(x,x)$, where $1 \leq l \leq r$. We call them \textit{unmixed circuits}.  
\begin{remark} \label{obs: obvio pero crucial}
    The unmixed circuit of greatest length $g_r(x,x)$ of $\widetilde{X}$ is a closed $2$-form on $\g_A$ that is non-degenerate on $\widetilde{X}$.
\end{remark}
\indent Denote by $\alt^2(\widetilde{X}, \widetilde{Y})^*$ the subspace of $\alt^2 \g_A^*$ spanned by all circuits of $\widetilde{X}$ and $\widetilde{Y}$, of all admissible lengths. For the case $\widetilde{X} = \widetilde{Y}$, we use the notation $\alt^2(\widetilde{X}, \widetilde{X})^*$.   
\begin{proposition} \label{prop: circuits}
    The set of cohomology classes of all circuits of $\widetilde{X}$ and $\widetilde{Y}$, of all admissible lengths, is a basis of $\alt^2(\widetilde{X}, \widetilde{Y})^*$. 
\end{proposition}
\begin{proof}
    Every circuit is closed either by direct computation or by Lemma \ref{lemma: X, Y dos}. Every circuit is non-exact, as exact $2$-forms are divisible by $f^1$ by Lemma \ref{lemma: elementary facts}(i); in fact, for the same reason, none of their linear combinations are exact. Clearly, the set of their cohomology classes is linearly independent. The fact that it also spans $\alt^2(\widetilde{X}, \widetilde{Y})$ follows from Lemma \ref{lemma: X, Y dos}. 
\end{proof}  
\begin{corollary} \label{cor: second betti, basically}
    The dimension of $\alt^2(\widetilde{X}, \widetilde{Y})^*$ is $r+s+2m$ if $\widetilde{X} \neq \widetilde{Y}$, and $r$ if $\widetilde{X} = \widetilde{Y}$. 
\end{corollary}
\indent Recall the decomposition of the spaces of forms of a fixed degree on $\g_A$ by $W_0$ and double generalized subspaces of $A_0$, introduced in equations \eqref{eq: g_A es f1 y u} and \eqref{eq: double generalized eigenspaces}:  
\begin{gather*}
    \alt^{\ell} \g_A^* = \R f^1 \otimes \alt^{\ell - 1} \u^* \oplus \alt^{\ell} \u^*, \\  
    \alt^{\ell} \u^* = \bigoplus_{i_0 + i_1 + \cdots + i_q = \ell} \alt^{i_0, i_1, \ldots, i_q}, \quad \text{where} \quad \alt^{i_0, \ldots, i_q} := \alt^{i_0} W_0^* \otimes \alt^{i_1} W_1^* \otimes \cdots \otimes \alt^{i_t} W_q^*,
\end{gather*} 
\noindent for all $0 \leq \ell \leq \dim \g_A$. Also, decompose every double generalized eigenspace $W_i$ ($i>0$) as a sum 
\begin{align*}
    W_i = \widetilde{X_{i,1}} \oplus \cdots \oplus \widetilde{X_{i,p_i}}
\end{align*}
\noindent of double elementary canonical subspaces, as in equation \eqref{eq: a finer decomposition}.

\indent Let $U$ be the subset of cohomology classes of non-exact forms $c f^1 \wedge \mu $ with $c \in \R$ and $\mu \in \u^*$, and $V$ the subspace in cohomology induced from $W_0^* \otimes W^*$, where $W := W_1 \oplus \cdots \oplus W_q$. 
\begin{theorem} \label{thm: second-degree cohomology}
    \phantom{.}
\begin{enumerate} [\rm (i)]
    \item $H^2(\g_A) = U \oplus V \oplus \alt^2 W_0^* \oplus \alt^2 W_1^* \oplus \cdots \oplus \alt^2 W_q^*$. 
    \item If $v = 0$ then $V = 0$ and $H^2(\g_A) = U \oplus \alt^2 W_0^* \oplus  \alt^2 W_1^* \oplus \cdots \oplus \alt^2 W_q^*$. 
    \item If $0 \notin \spec(A_0)$ then $\alt^2 W_0^* = \{0\}$ and $H^2(\g_A) = U \oplus \alt^2 W_1^* \oplus \cdots \oplus \alt^2 W_q^*$.
    \item For all $1 \leq i \leq q$, $\alt^2 W_i^* = \bigoplus_{a \neq b} \alt^2(\widetilde{X_{i,a}}, \widetilde{X_{i,b}})^*$. 
    \item $[f^1 \wedge f^2] \in U$. If further $0 \notin \spec(A_0)$ then $[f^1 \wedge f^2]$ is a basis of $U$.  
\end{enumerate}
\end{theorem}
\begin{proof} \phantom{.} \\ 
    \noindent (i) Every $2$-form splits as a sum of elements in the subspaces $\R f^1 \otimes \u^*$ and $\alt^2 \u^*$: 
    each term is either divisible by $f^1$ or not. It is then clear why the subspace $U$ appears in the decomposition of $H^2(\g_A)$. Moreover, due to Corollary \ref{cor: splitting}, every closed $2$-form on $\u$ splits as a sum of $2$-forms over the subspaces $W_0$, $W_1, \ldots, W_q$ and possibly terms mixing a $1$-form on $W_0$ with a $1$-form in any other $W_i$ ($i > 0$). Therefore, $V$ appears in the decomposition of $H^2(\g_A)$. Each term in such a splitting of a closed $2$-form is certainly closed, and they are all non-exact by Lemma \ref{lemma: elementary facts}(i).  

    \noindent (ii) If $v = 0$ then $V = \{0\}$ due to Corollary \ref{cor: splitting}(ii). The decomposition of $H^2(\g_A)$ is then a consequence of (i). 

    \noindent (iii) If $0 \notin \spec(A_0)$ then we may assume $v = 0$ in equation \eqref{eq: the alternative}, and moreover $W_0 = \R f_2$. Then $V = \{0\}$ as well as $\alt^2 W_0^* = \{0\}$, and the decomposition of $H^2(\g_A)$ follows from (i).  
    
    \noindent (iv) It follows from Lemma \ref{lemma: X, Y dos} and Proposition \ref{prop: circuits}.
    
    \noindent (v) It follows from Lemma \ref{lemma: elementary facts}(v) that $[f^1 \wedge f^2] \in U$. The second statement follows from Corollary \ref{cor: para citar}. \qedhere
\end{proof}

\indent In principle, Corollary \ref{cor: second betti, basically} and Theorem \ref{thm: second-degree cohomology} allow us to compute the second Betti number $b_2(\g_A)$ of $\g_A$ in some cases, linking it to the Jordan form of $A$. However, it appears that not much is gained in putting on effort to produce a somewhat explicit formula, as in this generality, it seems to be noninformative.  

\indent Of more use in what follows is to consider only the case where $A$ has a very simple Jordan form, described by stating that $\spec(A_0) = \{\lambda, -\lambda \}$ (with $\lambda\neq 0$) and $W_{\lambda}= \widetilde{X}$ is a double elementary canonical subspace of $A_0$ of size $m$. Note that, in this case, $\widetilde{X}= \u_0$; moreover, there is no $v$ in the alternative given in equation \eqref{eq: the alternative}. Setting $2n = \dim \g_A$, so $2n = 2m + 2$, it is possible to describe $H^1(\g_A)$, $H^2(\g_A)$, $H^{2n-2}(\g_A)$, and $H^{2n-1}(\g_A)$ in simple terms. In that regard, using the basis for $\widetilde{X}$ as described in equation \eqref{eq: x tilde base}, denote by $\Gamma = x^1 \wedge \cdots \wedge x^{2m}$ the top form in $\u_0$, and define 
\begin{gather*}
    \Gamma_a := x^1 \wedge \cdots \wedge \widehat{x^a} \wedge \cdots \wedge x^{2m} \in \alt^{2m-1} \u_0^*, \\
    \Gamma_{b, c} := x^1 \wedge \cdots \wedge \widehat{x^b} \wedge \cdots \wedge \widehat{x^c} \wedge \cdots \wedge x^{2m}  \in \alt^{2m - 2} \u_0^*, 
\end{gather*}
\noindent where $b < c$. The symbol $\widehat{x^a}$ means that the $1$-form $x^a$ does not appear in the expression for $\Gamma_a$, and similarly for $\Gamma_{b,c}$. Also, set $\delta := f^1 \wedge f^2$. This is a similar notation to that used by the authors in \cite{AG}. Note that $\delta$ and $\Gamma$ are closed forms, a fact that can be established by direct computation. For all $1 \leq l \leq m$, the \textit{companion $h_l(x,x)$ of the circuit $g_l(x,x)$} is defined to be
\begin{align*}
    h_l(x,x) := \sum_{i = 1}^l (-1)^{l + i + 1} \delta \wedge \Gamma_{i, m + l + 1 - i} \in \alt^{2n-2} \g_A^*.
\end{align*}
\noindent It is clear that $h_l$ is closed for all $1 \leq l \leq m$ according to Lemma \ref{lemma: elementary facts}(ii), since each $h_l$ is divisible by $f^1$. It is our immediate intention to prove that they are non-exact. For this we need to recall some facts concerning the Poincaré duality, and refer to \cite[Chapter 1, Section 3, pp 27]{Fuks} for details. Recall that it holds precisely when $\g_A$ is unimodular, and amounts to the fact that
\begin{gather*}
    H^k(\g_A) \times H^{2n-k}(\g_A) \to H^{2n}(\g_A) \cong \R \\
    (\theta, \eta) \mapsto \theta \wedge \eta  
\end{gather*}
\noindent is a non-degenerate bilinear pairing for all $k$ that make sense. The isomorphism $H^{2n}(\g_A) \cong \R$ results from viewing $\theta \wedge \eta$ as proportional to the volume form $\delta \wedge \Gamma$ on $\g_A$, say $\theta \wedge \eta = C \delta \wedge \Gamma$, and assigning $\theta \wedge \eta \mapsto C$. 

\indent In the next result we incur in the slight abuse of language of calling a class in cohomology by a representative of said class. 
\begin{theorem} \label{thm: h1, h2, h2n-2, h2n-1}
    If $0 \notin \spec(A_0)$ then
\begin{gather*}
    H^1(\g_A) = \Span_{\mathbb{R}} \{f^1, f^2\}, \quad H^{2n-1}(\g_A) = \Span_{\mathbb{R}} \{f^2 \wedge \Gamma, f^1 \wedge \Gamma\},     
\end{gather*}
    \noindent If further $\spec(A_0) = \{\lambda, -\lambda \}$ (with $\lambda >0$) and $W_{\lambda} = \widetilde{X}$ is a double elementary canonical subspace of $A_0$ of size $m$, then
\begin{gather*}
    H^2(\g_A) = \Span_{\mathbb{R}} \{ \delta, g_1(x,x), \ldots, g_m(x,x) \}, \quad 
    H^{2n-2}(\g_A) = \Span_{\mathbb{R}} \{\Gamma, h_1(x,x), \ldots, h_m(x,x) \}.  
\end{gather*}    
\end{theorem}
\begin{proof}
    It was noted throughout this section that all forms referred to in the statement of the Theorem are closed, and so they induce well-defined cohomology classes. 
    
    The hypotheses on $A$ imply that the only closed $1$-forms on $\g_A$ are $f^1$ and $f^2$, thereby establishing that $H^1(\g_A) = \Span_{\mathbb{R}} \{f^1, f^2\}$. From here and by Poincaré duality, the description of $H^{2n-1}(\g_A)$ is verified simply by noting that
\begin{align*}
    f^1 \wedge \left( f^2 \wedge \Gamma \right) = \delta \wedge \Gamma, \quad f^2 \wedge \left( f^1 \wedge \Gamma \right) = - \delta \wedge \Gamma,
\end{align*}
    \noindent implying that neither $f^1 \wedge \Gamma$ nor $f^2 \wedge \Gamma$ are non-exact. 

    \indent Assume now that the extra hypotheses are valid. The description for $H^2(\g_A)$ follows from Theorem \ref{thm: second-degree cohomology}, since there are only unmixed circuits. As before, the description of $H^{2n-2}(\g_A)$ follows from Poincaré duality, noting that 
\begin{align*}
    g_l(x,x) \wedge h_l(x,x) &= \left( \sum_{i = 1}^l (-1)^{i + 1} x^i \wedge x^{ m + l + 1 - i } \right) \wedge \left( \sum_{j = 1}^l (-1)^{l + j + 1} \delta \wedge \Gamma_{j, m + l + 1 - j} \right) \\
    &= \sum_{i = 1}^l \sum_{j = 1}^l (-1)^{l + i + j} \delta \wedge x^i \wedge x^{ m + l + 1 - i } \wedge \Gamma_{j, m + l + 1 - j} \\
    &= \sum_{i = 1}^l (-1)^{l + i + j} (-1)^{m+1} \delta \wedge \left( x^i \wedge \Gamma_i \right) \wedge \left( x^{ m + l + 1 - i } \wedge \Gamma_{m + l + 1 - i} \right) \\
    &= \sum_{i = 1}^l (-1)^{l + i + j} (-1)^{m+1} (-1)^{i-1} (-1)^{m+l+1-i-1} \delta \wedge \Gamma \\
    &= \sum_{i = 1}^l (-1)^{2m + 2l} \delta \wedge \Gamma = l \; \delta \wedge \Gamma
\end{align*}
    \noindent for all $1 \leq l \leq m$, implying that $h_l$ is non-exact for all $1 \leq l \leq m$. 
\end{proof}
\begin{remark}  
    Poincaré duality can be used to establish a decomposition for $H^{2n-2}(\g_A)$ similar to the one in Theorem \ref{thm: second-degree cohomology} for $H^2(\g_A)$. 
\end{remark}
\indent The following facts are central in the proof of the main result of this article. For them to hold we assume that $\spec(A_0) = \{\lambda, -\lambda\}$ (with $\lambda\neq 0$) and $W_{\lambda}= \widetilde{X}$ is a double elementary canonical subspace of $A_0$ of size $m \geq 2$. We use the rest of the notation introduced right before Theorem \ref{thm: h1, h2, h2n-2, h2n-1}, and also set $\rho := g_m(x,x)$, which is a closed $2$-form on $\g_A$, non-degenerate on $\widetilde{X}$ as per Remark \ref{obs: obvio pero crucial}.  
\begin{proposition} \label{prop: a direct computation 1} 
    For $m \geq 2$, $\rho^{m-2} \wedge (x^1 \wedge x^{m+1}) = \pm \Gamma_{m,2m}$. 
\end{proposition}
\begin{proof}
    Note that $g_m(x,x) = \sum_{i = 1}^m (-1)^{i + 1} x^i \wedge x^{ 2m + 1 - i }$ has precisely $m$ terms, and that the computation of $\rho^{m-2}$ involves choosing every possible combination of $m-2$ of these and wedge-multiplying them. Each term in $\rho^{m-2}$ that is divisible by either $x^1$ or $x^{m+1}$ multiplies to zero in the computation of $\rho^{m-2} \wedge (x^1 \wedge x^{m+1})$. Moreover, there is precisely one term in $\rho^{m-2}$ that is not divisible by either $x^1$ nor $x^{m+1}$, and it arises from multiplying all terms in $g_m(x,x)$ but the first and the last, those being $x^1 \wedge x^{2m}$ and $(-1)^{m+1} x^m \wedge x^{m+1}$. It is then clear that
\begin{align*}
    \rho^{m-2} \wedge (x^1 \wedge x^{m+1}) = \left( \bigwedge_{i = 2}^{m-1} (-1)^{i + 1} x^i \wedge x^{ 2m + 1 - i } \right) \wedge (x^1 \wedge x^{m+1}) 
\end{align*}
    \noindent is $\Gamma_{m,2m}$ up to a sign. 
\end{proof} 
\indent A similar reasoning as the one in Proposition \ref{prop: a direct computation 1} establishes that if 
\begin{align} \label{eq: necesito esta omega}
    \omega := \delta + g_m(x,x)
\end{align}
\noindent then $\omega^{n-2} \wedge (x^1 \wedge x^{m+1}) = \pm \delta \wedge \Gamma_{m,2m}$, where $n = m + 1$ and $m \geq 2$.  
\begin{proposition} \label{prop: a direct computation 2}
    For $m \geq 2$, $\delta \wedge \Gamma_{ m, 2m }$ is exact; moreover, $\delta \wedge \Gamma_{ m, 2m } = - d ( f^2 \wedge \Gamma_{ m, 2m-1} )$.  
\end{proposition}
\begin{proof} 
    The exactness of $\delta \wedge \Gamma_{ m, 2m }$ follows from Theorem \ref{thm: h1, h2, h2n-2, h2n-1}, since
\begin{align*}
    \{ \Gamma, h_1(x,x), \ldots, h_m(x,x), \delta \wedge \Gamma_{ m, 2m } \} \subseteq \alt^{2n-2} \g_A^*
\end{align*}
    \noindent is a linearly independent subset. The equality $\delta \wedge \Gamma_{ m, 2m } = - d ( f^2 \wedge \Gamma_{ m, 2m-1} )$ follows from direct computation. It is based on the facts that $d f^2 = 0$ as per Proposition \ref{prop: laura 2} (insofar we can take $v = 0$) and Lemma \ref{lemma: elementary facts}(iv), and $d \Gamma_{m, 2m-1} = f^1 \wedge \Gamma_{m, 2m}$, which can be established with the next relations:
\begin{itemize}
    \item For all $1 \leq i \leq m-1$,  
\begin{align*}
    x^1 \wedge \cdots \wedge (d x^i) \wedge \cdots \wedge x^{m-1} \wedge x^{m+1} \wedge \cdots \wedge x^{2m-2} \wedge x^{2m} 
    &= (-1)^i \lambda f^1 \wedge \Gamma_{m, 2m-1}.
\end{align*}    
    \item For $i = m+1$, 
\begin{align*}
    x^1 \wedge \cdots \wedge x^{m-1} \wedge (d x^{m+1}) \wedge \cdots \wedge x^{2m-2} \wedge x^{2m} 
    &= (-1)^{m-1} \lambda f^1 \wedge \Gamma_{m, 2m-1}. 
\end{align*}       
    \item For all $m+2 \leq i \leq 2m-2$,  
\begin{align*}
    x^1 \wedge \cdots \wedge x^{m-1} \wedge x^{m+1} \wedge \cdots \wedge (dx^i) \wedge \cdots \wedge x^{2m-2} \wedge x^{2m} 
    &= (-1)^i\lambda f^1 \wedge \Gamma_{m, 2m-1}. 
\end{align*}         
    \item For $i = 2m$,
\begin{align*}
    x^1 \wedge \cdots \wedge x^{m-1} \wedge x^{m+1} \wedge \cdots \wedge x^{2m-2} \wedge (d x^{2m}) 
    &= - \lambda f^1 \wedge \Gamma_{m, 2m-1} + f^1 \wedge \Gamma_{m, 2m}.  
\end{align*} 
\end{itemize}
\noindent One then gets 
\begin{align*}
    d \Gamma_{m, 2m-1} = &\phantom{+} \; \; \sum_{i=1}^{m-1} (-1)^{i+1} (-1)^i \lambda f^1 \wedge \Gamma_{m, 2m-1} + (-1)^{m+1} (-1)^{m-1} \lambda f^1 \wedge \Gamma_{m,2m-1} +\\
    &+ \sum_{i=m+2}^{2m-2} (-1)^i (-1)^i \lambda f^1 \wedge \Gamma_{m,2m-1} + (-1)^{2m-1} \left( - \lambda f^1 \wedge \Gamma_{m, 2m-1} + f^1 \wedge \Gamma_{m, 2m} \right) \\
    = &- \sum_{i=1}^{m-1} \lambda f^1 \wedge \Gamma_{m, 2m-1} + \lambda f^1 \wedge \Gamma_{m,2m-1} + \sum_{i=m+2}^{2m-2} \lambda f^1 \wedge \Gamma_{m,2m-1} + \\
    &+ \lambda f^1 \wedge \Gamma_{m, 2m-1} + (-1)^{2m-1} f^1 \wedge \Gamma_{m, 2m} \\
    =& \; \; - f^1 \wedge \Gamma_{m, 2m}. \qedhere
\end{align*}
\end{proof}
\begin{corollary} \label{cor: no es lefschetz}
    For $m \geq 2$ and $\omega$ as in equation \eqref{eq: necesito esta omega}, $\omega^{n-2} \wedge (x^1 \wedge x^{m+1})$ is an exact form.
\end{corollary} 

\subsection{The space of symplectic forms.} \label{section: the space of symplectic forms}

Consider a matrix $A$ as in equation \eqref{eq: the alternative}, with $a = 0$ and possibly with $v \neq 0$. Refer to Section \ref{section: notation} for the notation. In particular, recall that $\u = W_0 \oplus W$ where $W := \oplus_{i=1}^q W_i$ is the direct sum of all the double generalized eigenvalues of $A_0$. Also, recall the decomposition $A_0 = M \oplus N$, where $M := A \vert_{W_0}$ and $N := A_0 \vert_W$, under the assumption that $0 \in \spec(A_0)$. When $v \neq 0$, $M$ is nonzero and there is a (non-abelian) nilpotent subalgebra of $\g_A$ isomorphic to $\g_M$. 

\indent The results of the previous sections allow for a description of the structure of a generic symplectic form on $\g_A$. 
\begin{proposition} \label{prop: splitting of symplectic forms}
    Let $\omega$ be a symplectic form on $\g_A$ as in Proposition \ref{prop: LW}. 
\begin{enumerate} [\rm (i)]
    \item If $v \neq 0$ then there exist $\sigma \in \alt^2 \g_M^*$ and $\rho \in \alt^2 W^*$ such that $\omega = \sigma + \rho$. Both $\sigma$ and $\rho$ are closed $2$-forms on $\g_A$ and non-degenerate on the respective subspaces $\g_M$ and $W$ where they are defined. 
    \item If $v = 0$ then there exists a $1$-form $\mu \in \g_A^*$ and $2$-forms $\omega_i \in \alt^2 W_i^*$, with $0 \leq i \leq q$, such that $\omega = f^1 \wedge \mu + \sum_{i=0}^q \omega_i$. Moreover, for all $0 \leq i \leq q$, $\omega_i$ is a closed and non-exact $2$-form on $\g_A$ that is also non-degenerate on the subspace where it is defined. 
    \item If $v = 0$ then, for all $1 \leq i \leq q$, $\omega_i$ is a linear combination of circuits. If further $0 \notin \spec(A_0)$ then $\mu$ is proportional to $f^2$ and $\omega_0 = 0$.
\end{enumerate}    
\end{proposition}
\begin{proof} 
    Each splitting follows from the various statements in Corollary \ref{cor: splitting} and Theorem \ref{thm: second-degree cohomology}. The non-degeneracy of each term in each decomposition follows from the non-degeneracy of $\omega$ on $\g_A$ and Corollary \ref{cor: splitting}. Statement (iii) follows from Theorem \ref{thm: second-degree cohomology}(iv) and (v). 
\end{proof}
\begin{remark}
    We can interpret $f^1 \wedge f^2$ as a circuit of length $1$. With the further assumption that $0 \notin \spec(A_0)$, the statement of Proposition \ref{prop: splitting of symplectic forms}(iii) can be simplified to the following: every symplectic form on $\g_A$ is a linear combination of circuits plus an exact term.   
\end{remark}
\indent Assume now that $\spec(A_0) = \{\lambda, -\lambda\}$ (with $\lambda\neq 0$). Let
\begin{align*}
    W = \widetilde{X_1} \oplus \cdots \oplus \widetilde{X_p}
\end{align*}
\noindent  be the decomposition of the only double generalized eigenspace $W$ of $A_0$, associated to the pair of nonzero eigenvalues $\{\lambda, - \lambda\}$, as in equation \eqref{eq: a finer decomposition}. Concerning Proposition \ref{prop: splitting of symplectic forms}, and barring inaccuracies, if $\omega \in \alt^2 \g_A^*$ is a symplectic form on $\g_A$ then, up to scaling,  
\begin{align} \label{eq: generic omega}
    \omega = f^1 \wedge f^2 + \sum_{a \neq b} \sum_l b_l(x_a, x_b) g_l(x_a, x_b) + d \eta,  
\end{align}
\noindent for some scalars $b_l(x_a, x_b)$ and a $1$-form $\eta$ on $\g_A$. It is straightforward to verify that the non-degeneracy of $\omega$ is equivalent to the fact that each $b_m(x_a, x_b)$ is nonzero, by computing the volume form $\omega^n$. It is a small miracle that, with an appropriate change of basis, $\omega$ can be put in a much simpler form. 
\begin{theorem} \label{thm: splitting of symplectic forms}
    Let $\omega$ be a symplectic form on $\g_A$. Under the hypothesis and notations of the  paragraph above, there exists a basis $\{ \widetilde{x}_l \mid 1 \leq l \leq \dim \widetilde{X_i} \}$ of each subspace $\widetilde{X_i}$ such that
\begin{align} \label{eq: generic omega but better}
    \omega = f^1 \wedge f^2 + \sum_{i=1}^p g_{ m_i } ( \widetilde{x}_i, \widetilde{x}_i)  + d \eta
\end{align}
    \noindent Here, $g_{ m_i } ( \widetilde{x}_i, \widetilde{x}_i)$ represents the unmixed circuit of $\widetilde{X_i}$ with the largest admissible size $m_i := \dim X_i$, written in the new basis of $\widetilde{X_i}$. That is, possibly after a change of basis, every symplectic form on $\g_A$ is a sum of unmixed circuits of maximum length plus an exact term.
\end{theorem}
\begin{proof}
    We present a detailed argument only for the cases $p = 1$ and $p = 2$, in which the only double generalized eigenspace of $A_0$ associated to the nonzero pair $\{ \lambda, - \lambda\}$ is decomposed as either
\begin{align*}
    W = \widetilde{X} \text{ (for $p = 1$)} \quad \text{or} \quad W = \widetilde{X} \oplus \widetilde{Y} \text{ (for $p = 2$)},
\end{align*}
    \noindent since it allows for a good glimpse of the argument that works in general without the cumbersome notation. At the end we briefly discuss how the proof carries over for $p = 3$, where
\begin{align*}
    W = \widetilde{X} \oplus \widetilde{Y} \oplus \widetilde{Z} \text{ (for $p = 3$)},
\end{align*}
    \noindent clearing up the way for larger $p$.

    \indent We begin with the case $p = 1$. Write
\begin{align*}
    \omega = f^1 \wedge f^2 + \sum_{l=1}^m b_l g_l(x,x) + d \eta,
\end{align*}
    \noindent where $b_1, \ldots, b_m \in \mathbb{R}$, with $b_m$ nonzero, and $\eta \in \alt^1 \g_A^\ast$. The condition imposed on $b_m$ is equivalent to the non-degeneracy of $\omega$. Recall that
\begin{align*}
    g_l(x,x) = \sum_{i = 1}^l (-1)^{i + 1} x^i \wedge x^{ m + l + 1 - i }, \quad 1 \leq l \leq m. 
\end{align*}    
    \noindent The key step in the proof is the following computation: 
\begin{align*}
    \sum_{l=1}^m b_l g_l = & \sum_{l=1}^m b_l \sum_{i = 1}^{l} (-1)^{i + 1} x^i \wedge x^{ m + l + 1 - i } = \sum_{i=1}^m \sum_{l=i}^m b_l (-1)^{i + 1} x^i \wedge x^{ m + l + 1 - i }\\ & = \sum_{i = 1}^m (-1)^{i+1} x^i \wedge \left( \sum_{j = 1}^{m + 1 - i} b_{j + i - 1} x^{m + j} \right) = \sum_{i = 1}^m (-1)^{i+1} x^i \wedge \left( \sum_{j=1}^k b_{j + m - k} x^{m+j} \right). 
\end{align*}
    \noindent Inspired by the definition of a circuit, set 
\begin{align*}
    \widetilde{x}^{2m+1-i} = 
    \begin{cases}
        x^{2m + 1 - i}, & m + 1 \leq i \leq 2m, \\ 
        \sum_{j = 1}^{m + 1 - i} b_{j + i - 1} x^{ m + j }, & \phantom{m+} \; \, 1 \leq i \leq m, 
    \end{cases}   
\end{align*}
    \noindent so as to have
\begin{align*}
    \omega = f^1 \wedge f^2 + g_m( \widetilde{x}, \widetilde{x}) + d \widetilde{\eta}. 
\end{align*}

\noindent Note that the set $\{f^1,f^2\}\cup \{\widetilde{x}^j \mid 1\leq j\leq 2m\}$ is a basis of $\g_A^*$ since $\omega$ is non-degenerate.

\noindent Define $\varphi: \g_A^* \to \g_A^*$ to be the linear map arising from the transformations
\begin{align*}
    \varphi(f^k) = f^k \text{ for $k = 1, 2$}, \quad \varphi(x^i) = \widetilde{x}^i \text{ for $1 \leq j \leq 2m$}. 
\end{align*}
    \noindent It is straightforward to check that
\begin{align*}
    d \widetilde{x}^i = d x^i \text{ for all $1 \leq i \leq 2m$},
\end{align*}
    \noindent and therefore $\varphi$ gives rise to a Lie algebra homomorphism as per Lemma \ref{lemma: equivalencia para morfismos}. Since $\varphi$ preserves $\omega$, itself a non-degenerate $2$-form, $\varphi$ must be bijective and therefore $\{f^1, f^2\} \cup \{ \widetilde{x}^i \mid 1 \leq i \leq 2m\}$ is the dual of some basis of $\g_A$.  

    \indent We now tackle the case $p = 2$. There is no loss of generality in assuming $r \geq s$. Write
\begin{align*}
    \omega = f^1 \wedge f^2 + \sum_{l=1}^r b_l g_l(x,x) + \sum_{k=1}^s c_k g_k(y,y) + \sum_{m=1}^s d_m g_m(x,y) + \sum_{n=1}^s e_n g_n(y,x) + d \eta, 
\end{align*}  
    \noindent every symbol meaning the obvious. The non-degeneracy of $\omega$ is equivalent to the facts that $b_r$, $c_s$, $d_s$, and $e_s$ are all nonzero. As before, note that
\begin{align*}
    \sum_{l=1}^r b_l g_l(x,x) + \sum_{m=1}^s d_m g_m(x,y) &=
    \sum_{l=1}^r b_l \sum_{i = 1}^{l} (-1)^{i + 1} x^i \wedge x^{ r + l + 1 - i } + \sum_{m=1}^s d_m \sum_{i = 1}^m (-1)^{i + 1} x^i \wedge y^{ s + m + 1 - i } \\
    &= \sum_{i = 1}^r \sum_{l = i}^r (-1)^{i + 1} b_l x^i \wedge x^{ r + l + 1 - i } + \sum_{i = 1}^s \sum_{m = i}^s d_m (-1)^{i + 1} x^i \wedge y^{ s + m + 1 - i } \\
    &= \sum_{i = 1}^r (-1)^{i + 1} x^i \wedge \left( \sum_{l = i}^r b_l x^{ r + l + 1 - i } \right) + \sum_{i = 1}^s (-1)^{i + 1} x^i \wedge \left( \sum_{m = i}^s d_m y^{ s + m + 1 - i } \right) \\
    &= \sum_{i = 1}^r (-1)^{i + 1} x^i \wedge \left( \sum_{j = 1}^{r+1-i} b_{j + i - 1} x^{ r + j } \right) + \sum_{i = 1}^s (-1)^{i + 1} x^i \wedge \left( \sum_{j = 1}^{s + 1 - i} d_{j + i - 1} y^{ s + j } \right) 
\end{align*}
    \noindent In a similar fashion,  
\begin{align*}
    \sum_{l=1}^s c_l g_l(y,y) + \sum_{m=1}^s e_m g_m(y,x) &= \sum_{i = 1}^s (-1)^{i + 1} y^i \wedge \left( \sum_{j = 1}^{s+1-i} c_{j + i - 1} y^{ s + j } \right) + \sum_{i = 1}^s (-1)^{i + 1} y^i \wedge \left( \sum_{j = 1}^{s + 1 - i} e_{j + i - 1} x^{ s + j } \right) 
\end{align*}    
    \noindent The proof goes on as in the previous case. Set 
\begin{gather*}
    \tilde{x}^{2r+1-i} = 
    \begin{cases}
        x^{2r+1-i} & r + 1 \leq i \leq 2r, \\ 
        \sum_{j = 1}^{r+1-i} b_{j + i - 1} x^{ r + j } +  \sum_{j = 1}^{s + 1 - i} d_{j + i - 1} y^{ s + j } & \phantom{s+} \; \, 1 \leq i \leq r,
    \end{cases} \\ 
    \tilde{y}^{2s+1-i} = 
    \begin{cases}
        y^{2s+1-i} & s + 1 \leq i \leq 2s, \\ 
        \sum_{j = 1}^{s+1-i} c_{j + i - 1} y^{ s + j } +  \sum_{j = 1}^{s + 1 - i} e_{j + i - 1} x^{ s + j } & \phantom{s+} \; \, 1 \leq i \leq s.
    \end{cases}
\end{gather*}   
\noindent Clearly,
\begin{align*}
    \omega = f^1 \wedge f^2 + g_r(\tilde{x}, \tilde{x}) + g_s(\tilde{y}, \tilde{y}) + d \tilde{\eta}, 
\end{align*}    
    \noindent as well as
\begin{align*}
    d \tilde{x}^i = dx^i \text{ for all $1 \leq i \leq r$ }, \quad d \tilde{y}^l = dy^i \text{ for all $1 \leq l \leq s$}. 
\end{align*}   
    \noindent The proof concludes exactly as before. 

    \indent Lastly, we ponder briefly about the case $p = 3$, the only intention being to make it clear how to proceed with the proof for larger $p$. Note that, in the proof for the case $p = 2$, we had to set apart the circuits carrying the labels $(x,x)$ and $(x,y)$ from the circuits carrying the labels $(y,x)$ and $(y,y)$, and performing a ``factoring out"\! of sorts: For the circuits $(x,x)$ and $(x,y)$, the common factor is $\sum_{i=1}^r (-1)^{i+1} x^i$; for the circuits $(y,x)$ and $(x,x)$, the common factor is $\sum_{i=1}^r (-1)^{i+1} y^i$. Note that a similar procedure would have worked if we had performed this factoring out for $(x,x)$ and $(y,x)$ on the one side, and for $(y,y)$ and $(x,y)$ on the other side (the only difference being the new bases). If, for bookkeeping purposes, we arrange the four labels in a matrix as
\begin{align*}
    \left[ \begin{array}{c c}    
            (x,x) & (x,y) \\
            (y,x) & (y,y) 
        \end{array} \right],
\end{align*}
    \noindent then it is readily understood that the two factoring-out processes are ``row-like"\! or ``column-like"\!, both equally valid. The case $p = 3$ can be established by an analogous factoring-out procedure that can also be easily described in a matrix arrangement: If we put
\begin{align*}
    \left[ \begin{array}{c c c}    
            (x,x) & (x,y) & (x,z) \\
            (y,x) & (y,y) & (y,z) \\
            (z,x) & (z,y) & (z,z)
        \end{array} \right],
\end{align*}
    \noindent then we can factor out row-wise or column-wise. Either way, the new bases are immediate, and the proof works just as before. It is now evident how to manage the proof for larger $p$.   
\end{proof}

\begin{remark}\label{rem: mas autovalores}
The role of the hypothesis on the spectrum of $A$ in Theorem \ref{thm: splitting of symplectic forms} is only to clarify its statement and its proof. Concisely, what was proven there amounts to the fact that every symplectic form on a unimodular almost abelian Lie algebra with $0 \notin \spec(A_0)$, possibly after a change of basis, can be written as a sum of essentially two factors: one containing multiple terms, one for each \textit{elementary} Jordan block associated to a nonzero eigenvalue of $A$, equal to an unmixed circuit of maximum length; other, being exact and therefore unimportant for cohomological purposes. 
\end{remark}

\begin{remark} \label{obs: cohomological uniqueness?}
    It is instructive to compare Theorem \ref{thm: splitting of symplectic forms} with \cite[Theorem 1.1]{CM}, where it is shown that diagonal almost abelian Lie algebras $\g_A$ admit at most one symplectic structure up to automorphism and scaling, and that a necessary (but not sufficient) condition for them to exist is that $\g_A$ be unimodular. Moreover, as per Theorem \ref{thm: second-degree cohomology}, it is possible to interpret that result as saying that any symplectic form on a diagonal unimodular almost abelian Lie algebra must be equivalent to one splitting as a sum of (necessarily unmixed) circuits of length $1$, one of them being precisely $f^1 \wedge f^2$. 
\end{remark}

\subsection{The failure of the hard-Lefschetz condition.} 

\indent Consider a matrix $A$ as in equation \eqref{eq: the alternative}, with  $a = 0$ and possibly with $v \neq 0$. In particular, recall that $\u = W_0 \oplus W$ where $W := \oplus_{i=1}^q W_i$ is the direct sum of all the double generalized eigenvalues of $A_0$. Also, recall the decomposition $A_0 = M \oplus N$, where $M := A \vert_{W_0}$ and $N := A_0 \vert_W$, under the assumption that $0 \in \spec(A_0)$. When $v \neq 0$, $M$ is nonzero and there is a (non-abelian) nilpotent subalgebra of $\g_A$ isomorphic to $\g_M$.
 
Although the main result of this section ultimately does not depend on the alternative in equation \eqref{eq: the alternative}, its proof does. In the case $0 \in \spec(A_0)$, where $A = M \oplus N$, $M$ is a nilpotent matrix, and therefore there is a nilpotent subalgebra of $\g_A$ isomorphic to $\g_M$. 
\begin{theorem} \label{thm: main result 1}
     Let $\g_A$ be unimodular, and assume additionally that $0 \in \spec(A_0)$. If $M$ is not the zero matrix (in particular, if $v \neq 0$ in equation \eqref{eq: the alternative}) then any symplectic form $\omega$ on $\g_A$ fails to satisfy the hard-Lefschetz condition. Moreover, the failure occurs at degree $1$. 
\end{theorem} 
\begin{proof}
    We established in the last paragraph that there is a nilpotent subalgebra $\g_M$ within $\g_A$, and that $W$ is an ideal in $\g_M$. Thus, it makes sense to express $\g_A = \g_M \ltimes W$. Note that $\g_M$ is non-abelian since $M$ is assumed to be a non-zero matrix. According to Proposition \ref{prop: splitting of symplectic forms}(i), any symplectic form $\omega$ on $\g_A$ can be decomposed as $\omega = \sigma + \rho$, where $\sigma \in \alt^2 \g_M^*$ and $\rho \in \alt^2 W^*$, and in particular $\g_M$ and $W$ are $\omega$-orthogonal. Moreover, $\sigma$ is a symplectic form on $\g_M$ and $\rho$ is a closed non-degenerate $2$-form on $W$. Consequently, Theorem \ref{thm: Benson y Gordon} applies, and so $\sigma$ does not satisfy the hard-Lefschetz condition on $\g_M$; moreover, the same result ensures that failure for this condition happens at degree $1$. By propagation as discussed in Lemma \ref{lemma: propagacion algebraica}, the same conclusion holds for $\omega = \sigma + \rho$ on $\g_A$.    
\end{proof}
\indent The case where $0 \notin \spec(A_0)$ and $M$ is the zero matrix is dealt with either by Theorem \ref{thm: Kasuya} or by the next result, depending on whether $N$ has nilpotent part or not. 
\begin{theorem} \label{thm: main result 2}
     Let $\g_A$ be unimodular and assume additionally that $0 \notin \spec(A_0)$. If $A_0$ has nonzero nilpotent part then any symplectic form $\omega$ on $\g_A$ fails to have the hard-Lefschetz condition. Moreover, the failure occurs at degree $2$, and never at degree $1$.  
\end{theorem} 
\begin{proof} 
    Let $\Xi = \{ \lambda_l \mid 1 \leq l \leq q \}$ be the set of positive eigenvalues of $A_0$. As in Section \ref{section: notation}, consider the decompositions
\begin{align*}
    \g_A = \Span_{\R} \{f_1, f_2\} \oplus W_1 \oplus \cdots \oplus W_q, \quad W_i = \widetilde{X_{i,1}} \oplus \cdots \oplus \widetilde{X_{i,p_i}} \text{ for all $1 \leq i \leq q$}. 
\end{align*}
    \noindent Note that we have already used the fact that $0 \notin \spec(A_0)$. Notice that the hypotheses on $A_0$ imply that it has a positive eigenvalue $\lambda_l > 0$ for which there is a Jordan elementary canonical block $J_m(\lambda_l)$ with $m \geq 2$. This Jordan elementary canonical block is associated with a double elementary canonical subspace $\widetilde{X}_{l,\ell}$ appearing in the decomposition of $W_l$. Set 
\begin{align*}
    \h_1 := \Span_{\R}\{f_1,f_2\} \oplus \widetilde{X_{l,\ell}}, \quad \h_2 := \bigoplus_{i \neq l} W_i \oplus \bigoplus_{j \neq \ell} \widetilde{X_{l,j}}.
\end{align*}
    \noindent Notice that $\h_1$ is a subalgebra of $\g_A$ and that $\h_2$ is an abelian ideal in $\g_A$, and so $\g_A = \h_1 \ltimes \h_2$. As indicated in Remark \ref{rem: mas autovalores}, setting $\omega_1 := f^1 \wedge f^2 + g_m( \widetilde{x}_{l, \ell}, \widetilde{x}_{l, \ell} )$, every symplectic form $\omega$ on $\g_A$ can be written as
\begin{align} \label{eq: are we done yet}
    \omega = \omega_1 + \omega_2 + d \eta = f^1 \wedge f^2 + g_m( \widetilde{x_{l, \ell}}, \widetilde{x_{l, \ell}} ) + \omega_2 + d \eta, 
\end{align}
    \noindent where $\omega_2$ is a linear combination of unmixed circuits of maximal lengths associated with the double elementary canonical subspaces within $\h_2$, and $\eta$ can be further assumed to be zero since we intend to prove a cohomological statement. This shows that $\omega$ splits as a symplectic form $\omega_1 = f^1 \wedge f^2 + g_m( \widetilde{x_{l, \ell}}, \widetilde{x_{l, \ell}} )$ on $\h_1$, and a closed non-degenerate $2$-form $\omega_2$ on $\h_2$. Thus, propagation in the sense of Lemma \ref{lemma: propagacion algebraica} applies, and therefore the failure of $(\g_A, \omega)$ to satisfy the hard-Lefschetz condition at degree $2$ can be established simply by arguing that $(\h_1, \omega_1)$ does not satisfy the hard-Lefschetz condition at degree $2$. As $2m = \dim \widetilde{X_{l, \ell}}$ with $m \geq 2$, Corollary \ref{cor: no es lefschetz} applies. This shows that the Lefschetz operator $L_{m-1}:H^2(\h_1) \to H^{2m}(\h_1)$ at degree $2$, given by $L_{m-1}(\cdot) := \frac{1}{(m-1)!} [\omega_1^{m - 1} \wedge \cdot \;]$, has a nontrivial kernel. This establishes the first part of the proof.

    It remains to see that the Lefschetz operator $L_{n-1}:H^1(\g_A) \to H^{2n-1}(\g_A)$ at degree $1$ is a bijection for a symplectic form as in equation \eqref{eq: are we done yet}, although now it is more convenient to write it as
\begin{align*}
    \omega = f^1 \wedge f^2 + \omega_0 + d \eta,
\end{align*}
    \noindent where $\omega_0$ is a linear combination of unmixed circuits of maximal length. As before, we take $\eta = 0$. Recall the descriptions of $H^1(\g_A)$ and $H^{2n-1}(\g_A)$ as in Theorem \ref{thm: h1, h2, h2n-2, h2n-1}. Since circuits of maximal length are non-degenerate on the subspaces where they are defined, it is clear that 
\begin{align*}
    \omega^{n-1} \wedge f^1 = \omega_0^{n-1} \wedge f^1 = \Gamma \wedge f^1, \quad \omega^{n-1} \wedge f^2 = \omega_0^{n-1} \wedge f^2 = \Gamma \wedge f^2,
\end{align*}    
    \noindent and therefore $L_{n-1}$ takes $H^1(\g_A)$ to $H^{2n-1}(\g_A)$ bijectively. 
\end{proof}

\begin{remark}
    Completely solvable almost abelian Lie algebras $\g_A$ as those described in Theorem \ref{thm: main result 2} are $1$-Lefschetz, by which we mean the Lefschetz operator $L_{n-1}:H^1(\g_A) \to H^{2n-1}(\g_A)$ at degree $1$ is a bijection for all completely solvable almost abelian Lie algebras $\g_A$. All these Lie algebras are shown to be of the form $\R^2 \ltimes \R^{2m}$ for some $m \in \N$, exhibiting the properties listed in \cite[Theorem 2]{BG2}.
\end{remark}

\indent Together with Theorem \ref{thm: Kasuya}, Theorems \ref{thm: main result 1} and \ref{thm: main result 2} establish the main result of this article.
\begin{theorem} \label{thm: principal}
    Let $\omega$ be a symplectic form on a unimodular completely solvable almost abelian Lie algebra $\g_A$. Then $\omega$ satisfies the hard-Lefschetz condition if and only if $A$ is semisimple. 
\end{theorem}

    
\section{Existence of lattices} \label{Section: lattices}

\indent In this section we use Proposition \ref{prop: useful criterion for lattices} to establish the existence of lattices in the simply connected Lie groups $G_A$ associated with some symplectic completely solvable almost abelian Lie algebras $\g_A$. 

\indent Consider a matrix $A$ as in equation \eqref{eq: the alternative} with $a = 0$ and $v = 0$. As before, we set $\dim \g_A = 2n$ and $m := n - 1$; in particular, $A_0 \in \mathfrak{sp}(m, \R)$. Proposition \ref{prop: LW} ensures that $\g_A$ is unimodular and that it admits a symplectic form.  We use Proposition \ref{prop: laura 2} in order to specify the symplectic matrices in which we are interested. The real numbers
\begin{align}\label{eq: tm}
    t_k:=\log \frac{k+\sqrt{k^2-4}}{2}, \quad \text{$k \in \N$ with $k \geq 3$}, 
\end{align}
\noindent play a fundamental role in the ensuing construction. Note that $\exp(t_k)$ and $\exp(-t_k)$ are the only roots of $p_k(x) := x^2 - kx + 1 \in \Z[x]$. We divide the construction in cases for a cleaner treatment. 

\medskip 

\textsl{Case (i):} Let $t \geq 1$ be an integer. Assume that $\spec(A_0) = \{0\}$ and that
\begin{align*}
    A_0 = J_{2t}(0) \in \mathfrak{sp}(t, \R).
\end{align*}
\noindent That $A_0$ belongs to $\mathfrak{sp}(t, \R)$ follows from Proposition \ref{prop: laura 2}. Since $\g_A$ is a nilpotent algebra with rational structure constants (in fact, only $0$’s and $1$’s), Malcev’s criterion (Theorem \ref{thm: Malcev}) guarantees the existence of a lattice in $\g_A$. However, this approach is overkill and does not fit well within our argument, so we give a more elementary proof: Since $\exp(A) = (1) \oplus \exp(A_0)$ is a unipotent matrix such that $\exp(A) - \mathrm{Id}$ is nilpotent of rank $2t$, both its characteristic and minimal polynomials coincide, and are given by
\begin{align*}
    (x - 1)^{2t+1} \in \Z[x].
\end{align*}
\noindent Consequently, $\exp(A)$ is conjugate to the companion matrix of $(x - 1)^{2t+1} \in \Z[x]$, which is an integer matrix. By Proposition \ref{prop: useful criterion for lattices}, the corresponding Lie group $G_A$ admits a lattice $\Gamma$. Note that the symplectic solvmanifold $\Gamma \backslash G_A$ does not satisfy the hard-Lefschetz condition, as follows from Theorem \ref{thm: main result 1}. (in actuality, it follows directly from Benson and Gordon's Theorem \ref{thm: Benson y Gordon}). 

\textsl{Case (ii):} Let $k \geq 3$ and $m \geq 2$ be integers. Assume that $\spec(A_0) = \{t_k, -t_k\}$, where $t_k$ is given by \eqref{eq: tm}, and that
\begin{equation*}
    A_0 = J_m(t_k) \oplus J_m(-t_k) \in \mathfrak{sp}(m,\R).
\end{equation*}
\noindent As before, the fact that $A_0 \in \mathfrak{sp}(m, \R)$ follows from Proposition \ref{prop: laura 2}. Moreover, since $\exp(A) = (1) \oplus \exp(A_0)$, both the characteristic and minimal polynomials of $\exp(A_0)$ coincide and are given by
\begin{equation*}
    (x-\e^{t_k})^m (x-\e^{-t_k})^m = (p_k(x))^m \in \Z[x],
\end{equation*}
\noindent where $p_k(x) = x^2 - kx + 1$. By a standard fact from linear algebra, $\exp(A_0)$ is conjugate to the companion matrix of $(p_k(x))^m \in \Z[x]$, which is an integer matrix. Thus, $\exp(A)$ is conjugate to an integer matrix, and by Proposition \ref{prop: useful criterion for lattices}, the corresponding Lie group $G_A$ admits a lattice $\Gamma$. The condition $m \geq 2$ ensures that the symplectic solvmanifold $\Gamma \backslash G_A$ does not satisfy the hard-Lefschetz condition, as follows from Theorem \ref{thm: main result 2}.

\textsl{Case (iii):} Let $t\geq 0$, $k_1, \dots, k_r\geq 3$, and $m_1, \dots, m_r \geq 1$ be integers satisfying the additional conditions that the $k_j$'s are distinct and that either $t \neq 0$ or $m_i \geq 2$ for some $1 \leq i \leq r$. Define $m := t + m_1 + \cdots + m_r$ and assume that
\begin{align*}
    A_0 = J_{2t}(0) \oplus J_{m_1}(t_{k_1}) \oplus J_{m_1}(-t_{k_1}) \oplus \cdots \oplus J_{m_r}(t_{k_r}) \oplus J_{m_r}(-t_{k_r}) \in \mathfrak{sp}(m,\R).
\end{align*}
\noindent Again, the fact that $A_0 \in \mathfrak{sp}(m, \R)$ follows from Proposition \ref{prop: laura 2}. Similarly to the previous cases,
\begin{align*}
    \exp(A)=(1) \oplus \exp[ J_{2t}(0) ]  \oplus \exp[J_{m_1}(t_{k_1})\oplus J_{m_1}(-t_{k_1})]\oplus \cdots \oplus \exp[J_{m_r}(t_{k_r})\oplus J_{m_r}(-t_{k_r})]. 
\end{align*}
\noindent The block $\exp[J_{2t}(0)]$ is conjugate to an integer matrix by Case (i), while each block $\exp[J_{m_i}(t_{k_i}) \oplus J_{m_i}(-t_{k_i})]$ is conjugate to an integer matrix by Case (ii). Therefore, $\exp(A)$ itself is conjugate to an integer matrix, implying that the corresponding Lie group $G_A$ admits a lattice $\Gamma$. Since either $t \neq 0$ or $m_i \geq 2$ for some $1 \leq i \leq r$, the symplectic solvmanifold $\Gamma \backslash G_A$ does not satisfy the hard-Lefschetz condition, as a consequence of either Theorem \ref{thm: main result 1} or Theorem \ref{thm: main result 2}.

\ 



\textit{Acknowledgements.} This work was partially supported by CONICET, SECyT-UNC and ANPCyT (Argentina). The authors are grateful to Romina M. Arroyo,  María Laura Barberis and Hisashi Kasuya for useful comments and suggestions.

\


\begin{thebibliography}{99}

\bibitem{AG}
A.\ Andrada and A.\ Garrone. Construction of symplectic solvmanifolds satisfying the hard-Lefschetz condition. \textit{Linear Algebra Appl.} \textbf{706} (2025), 70--100. 

\bibitem{AK}
D.\ Angella and H.\ Kasuya. Symplectic Bott-Chern cohomology of solvmanifolds. \textit{J. Symplectic Geom.} \textbf{17} (2019), 41--91. 

\bibitem{etal-nilp}
R.\ M.\ Arroyo, M.\ L.\ Barberis, V.\ Díaz, Y.\ Godoy and I.\ Hernández. Classification of nilpotent almost abelian Lie groups admitting left-invariant complex or symplectic structures. Preprint 2024, arXiv:2406.06819.

\bibitem{etal}
R.\ M.\ Arroyo, M.\ L.\ Barberis, V.\ Díaz, Y.\ Godoy and I.\ Hernández. Classification of almost abelian Lie groups admitting left-invariant complex or symplectic structures.  In preparation.

\bibitem{BT}
P.~de Bartolomeis and A.~Tomassini. On solvable generalized Calabi-Yau manifolds. \textit{Ann. Inst. Fourier} \textbf{56} (2006), 1281--1296.

\bibitem{BG1}
C.~Benson and C.~S.~Gordon. Kähler and symplectic structures on nilmanifolds. \textit{Topology} \textbf{27} (1988), 513--518. 

\bibitem{BG2}
C.~Benson and C.~S.~Gordon. Kähler structures on compact solvmanifolds. \textit{Proc. Am. Math. Soc.} \textbf{108} (1990), 971--980.

\bibitem{Bock} 
C.\ Bock. On low-dimensional solvmanifolds. \textit{Asian J. Math.} \textbf{20} (2016), 199--262. 

\bibitem{Brylinski} 
J.~L.~Brylinski. A differential complex for Poisson manifolds. \textit{J. Differ. Geom.} \textbf{28} (1988), 93--114. 

\bibitem{CM}
L.~P.~Castellanos Moscoso. Left-invariant symplectic structures on diagonal almost abelian Lie groups. \textit{Hiroshima Math. J.} \textbf{52} (2022), 357--378. 

\bibitem{Cavalcanti 1}
G.\ R.\ Cavalcanti. New aspects of the $dd^c$-Lemma. Oxford Thesis D. Phil. (2005), arxiv.org/abs/math/0501406.

\bibitem{Cavalcanti 2}
G.~R.~Cavalcanti. The Lefschetz property, formality, and blowing up in symplectic geometry. \textit{Trans. Am. Math. Soc.}  \textbf{359} (2007), 333--348.

\bibitem{Cho}
Y.~Cho. Hard-Lefschetz condition of symplectic structures on compact Kähler manifolds. \textit{Trans. Am. Math. Soc.} \textbf{368} (2016), 8223--8248. 

\bibitem{Freibert}
M.~Freibert. Cocalibrated structures on Lie algebras with a codimension-one abelian ideal. \textit{Ann. Glob. Anal. Geom.} \textbf{42} (2012), 537--563. 

\bibitem{Fuks}
D.~B.~Fuks. \textit{Cohomology of infinite-dimensional Lie algebras}. Springer, Monographs in Contemporary Mathematics (1986). 


\bibitem{Hattori}
A.~Hattori. Spectral sequence in the de Rham cohomology of fibre bundles. \textit{J. Fac. Sci. Univ. Tokyo Sect. I} \textbf{8} (1960), 289--331. 

\bibitem{Kasuya}
H.~Kasuya. Formality and hard-Lefschetz condition of aspherical manifolds. \textit{Osaka J. Math.} \textbf{50} (2013), 439--455. 

\bibitem{Kasuya-JDG}
H.~Kasuya. Minimal models, formality, and hard Lefschetz properties of solvmanifolds with local systems. \textit{J. Differ. Geom.} \textbf{93} (2013), 269--297.

\bibitem{LRV} 
J.\ Lauret and E.\ Rodr\'iguez-Valencia. On the Chern-Ricci flow and its solitons for Lie groups, \textit{Math. Nachr.} \textbf{288} (2015), 1512--1526. 

\bibitem{LW} 
J.\ Lauret and C.\ Will. On the symplectic curvature flow for locally homogeneous manifolds. \textit{J. Symplectic Geom.} \textbf{15} (2014), 1--49. 

\bibitem{LT}
F.\ Lusetti and A.\ Tomassini. Hard Lefschetz condition on symplectic non-Kähler solvmanifolds. Preprint 2025, arXiv:2501.13179. 

\bibitem{Malcev}
A.~Malcev. On a class of homogeneous spaces. \textit{Izv. Akad. Nauk SSSR} \textbf{13} (1949), 9--32; English translation in \textit{Am. Math. Soc. Transl.} \textbf{39} (1951). 

\bibitem{Mathieu}
O.~Mathieu. Harmonic cohomology classes of symplectic manifolds. \textit{Comment. Math. Helv.} \textbf{70} (1995), 1--9.


\bibitem{MMRR}
C.\ Mehl, V.\ Mehrmann, A.\ C.\ M.\ Ran and L.\ Rodman. Eigenvalue perturbation theory of structured real matrices and their sign characteristics under generic structured rank-one perturbations, \textit{Linear Multilinear Algebra} \textbf{64} (2016), 527--556.

\bibitem{Merkulov}
S.\ A.\ Merkulov. Formality of canonical symplectic complexes and Frobenius manifolds. \textit{Internat. Math. Res. Notices.} \textbf{1998} (1998), 727--733.

\bibitem{Milnor} 
J.~Milnor. Curvatures of left invariant metrics on Lie groups. \textit{Adv. Math.} \textbf{21} (1976), 293--329. 

\bibitem{Milovanov}
M.~V.~Milovanov, A description of solvable Lie groups with a given uniform subgroup, \textit{Math. USSR Sbornik} \textbf{41} (1982), 83--99.

\bibitem{Mostow}
G.~D.~Mostow. Factor spaces of solvable groups. \textit{Ann. Math.} \textbf{60} (1954), 1--27. 

\bibitem{Nomizu}
K.~Nomizu. On the cohomology of compact homogeneous spaces of nilpotent Lie groups. \textit{Ann. Math.} \textbf{59} (1954), 531--538. 


\bibitem{Saito}
M.~Saito. Sur certains groupes de Lie résolubles II. \textit{Sci. Papers Coll. Gen. Ed. Univ. Tokyo} \textbf{7} (1957), 157--168. 

\bibitem{Sawai}
H.~Sawai. A construction of lattices on certain solvable Lie groups, \textit{Topology Appl.} \textbf{154} (2007), 3125--3134.

\bibitem{SY}
H.~Sawai and T.~Yamada. Lattices on Benson–Gordon type solvable Lie groups. \textit{Topology Appl}. \textbf{149} (2005), 85–95

\bibitem{TT} 
Q.~Tan and A.~Tomassini. Remarks on some compact symplectic solvmanifolds. \textit{Acta. Math. Sin., English Ser.} \textbf{39} (2023), 1874--1886. 

\bibitem{Witte} 
D.~Witte. Superrigidity of lattices in solvable Lie groups, \textit{Invent. Math.} \textbf {122} (1995), 147--193.

\bibitem{Yamada 1}
T.~Yamada. Examples of compact Lesfschetz manifolds. \textit{Tokyo J. Math.} \textbf{25} (2002), 261--283.

\bibitem{Yamada 2}
T.~Yamada. A construction of lattices in splittable solvable Lie groups. \textit{Kodai Math. J.} \textbf{39} (2016), 378--388.

\bibitem{Yan}
D.~Yan. Hodge structure on symplectic manifolds. \textit{Adv. Math.} \textbf{120} (1996), 143–154.


\end{thebibliography}
\end{document}